\documentclass[11pt,a4paper]{article}
\usepackage{amsmath}
\usepackage{amssymb}
\usepackage{array}
\usepackage{enumerate}
\usepackage{algorithm}
\usepackage{algpseudocode}
\usepackage{float}

\usepackage{bm}
\usepackage[font=small,labelfont=bf]{caption}
\usepackage{subcaption}
\usepackage{changepage}
\usepackage[percent]{overpic}
\usepackage{mathtools}
\usepackage{amsfonts}
\usepackage{graphicx}
\usepackage{stackrel}
\usepackage{epstopdf}
\ifpdf
  \DeclareGraphicsExtensions{.eps,.pdf,.png,.jpg}
\else
  \DeclareGraphicsExtensions{.eps}
\fi
\usepackage{amsmath}
\usepackage[amsmath,thmmarks,hyperref]{ntheorem}
\usepackage{xr-hyper}
\usepackage{ifpdf}
\usepackage{hyperref}
\usepackage[capitalize,nameinlink]{cleveref}
\usepackage[all]{hypcap}
\usepackage{graphics,graphicx}
\usepackage{xcolor}

\theoremstyle{plain}
\theoremheaderfont{\normalfont\sc}
\theorembodyfont{\normalfont\itshape}
\theoremseparator{.}
\theoremsymbol{}

\newcommand{\proofbox}{\vbox{\hrule height0.6pt\hbox{\vrule height1.3ex width0.6pt\hskip0.8ex\vrule width0.6pt}\hrule height0.6pt}}

\theoremstyle{nonumberplain}
\theoremheaderfont{\normalfont\itshape}
\theorembodyfont{\normalfont}
\theoremseparator{.}
\theoremsymbol{\proofbox}

\newtheorem{thrm}{Theorem}[section]

\theoremstyle{definition}

\newtheorem{rmrk}[thrm]{Remark}

\theoremstyle{nonumberplain}
\theoremheaderfont{\normalfont\itshape}
\theorembodyfont{\normalfont}
\theoremseparator{.}
\theoremsymbol{\proofbox}

\newcommand{\newremark}[2]{
  \theoremstyle{plain}
  \theoremheaderfont{\normalfont\itshape}
  \theorembodyfont{\normalfont}
  \theoremseparator{.}
  \theoremsymbol{}
  \newtheorem{#1}[theorem]{#2}
}

\newremark{assumption}{Assumption}
\newremark{example}{Example}
\crefname{assumption}{Assumption}{Assumption}
\creflabelformat{enumi}{(#2#1#3)}
\crefname{enumi}{}{}

\usepackage{amsopn}
\newcommand{\xR}{\mathbb R}
\newcommand{\xC}{\mathbb C}
\newcommand{\xZ}{\mathbb Z}
\newcommand{\xN}{\mathbb N}
\def\xLtwo{{\rm L}^{2}}
\def\xHone{{\rm H}^{1}}

\newcommand{\nn}{\mathbf n}
\newcommand{\pp}{\mathbf p}
\newcommand{\ssigma}{\bm{\sigma}}
\newcommand{\ppsi}{\bm{\psi}}
\newcommand{\PP}{\bm{\mathcal{P}}}
\newcommand{\abs}[1]{\left|#1\right|}
\newcommand{\norm}[1]{\left\|#1\right\|}
\DeclareMathOperator*{\argmin}{arg\,min}

\title{Non-intrusive double-greedy parametric model reduction by interpolation of frequency-domain rational surrogates}
\author{Fabio Nobile and Davide Pradovera\thanks{Corresponding author. The authors acknowledge support from Swiss National Science Foundation project 182236.}\\[.25cm] CSQI, EPFL, Lausanne, Switzerland\\ \texttt{\{fabio.nobile,davide.pradovera\}@epfl.ch}}

\crefname{section}{section}{sections}
\crefname{subsection}{subsection}{subsections}
\Crefname{section}{Section}{Sections}
\Crefname{subsection}{Subsection}{Subsections}
\Crefname{figure}{Figure}{Figures}
\crefformat{equation}{\textup{#2(#1)#3}}
\crefrangeformat{equation}{\textup{#3(#1)#4--#5(#2)#6}}
\crefmultiformat{equation}{\textup{#2(#1)#3}}{ and \textup{#2(#1)#3}}
{, \textup{#2(#1)#3}}{, and \textup{#2(#1)#3}}
\crefrangemultiformat{equation}{\textup{#3(#1)#4--#5(#2)#6}}%
{ and \textup{#3(#1)#4--#5(#2)#6}}{, \textup{#3(#1)#4--#5(#2)#6}}{, and \textup{#3(#1)#4--#5(#2)#6}}
\Crefformat{equation}{#2Equation~\textup{(#1)}#3}
\Crefrangeformat{equation}{Equations~\textup{#3(#1)#4--#5(#2)#6}}
\Crefmultiformat{equation}{Equations~\textup{#2(#1)#3}}{ and \textup{#2(#1)#3}}
{, \textup{#2(#1)#3}}{, and \textup{#2(#1)#3}}
\Crefrangemultiformat{equation}{Equations~\textup{#3(#1)#4--#5(#2)#6}}%
{ and \textup{#3(#1)#4--#5(#2)#6}}{, \textup{#3(#1)#4--#5(#2)#6}}{, and \textup{#3(#1)#4--#5(#2)#6}}
\crefdefaultlabelformat{#2\textup{#1}#3}

\usepackage[margin = 1.in]{geometry}
\numberwithin{equation}{section}

\usepackage{tikz}
\usepackage[siunitx]{circuitikz}
\usepackage{pgfplots}
\usetikzlibrary{plotmarks, shapes.geometric}
\pgfplotsset{compat=newest}
\pgfplotsset{tick label style={font=\footnotesize}}
\pgfplotsset{colormap={whiteblack}{color(0)=(white) color(1)=(black)}}
\usepackage{url}
\usepackage{bm}
\usepackage{enumerate}
\usepackage{algorithm}
\usepackage{algpseudocode}

\ifpdf
\hypersetup{
  pdftitle={Non-intrusive double-greedy parametric model reduction by interpolation of frequency-domain rational surrogates},
  pdfauthor={Davide Pradovera}
}
\fi

\begin{document}
\maketitle

\begin{abstract}
We propose a model order reduction approach for non-intrusive surrogate modeling of parametric dynamical systems. The reduced model over the whole parameter space is built by combining surrogates in frequency only, built at few selected values of the parameters. This, in particular, requires matching the respective poles by solving an optimization problem. If the frequency surrogates are constructed by a suitable rational interpolation strategy, frequency and parameters can both be sampled in an adaptive fashion. This, in general, yields frequency surrogates with different numbers of poles, a situation addressed by our proposed algorithm. Moreover, we explain how our method can be applied even in high-dimensional settings, by employing locally-refined sparse grids in parameter space to weaken the curse of dimensionality. Numerical examples are used to showcase the effectiveness of the method, and to highlight some of its limitations in dealing with unbalanced pole matching, as well as with a large number of parameters.
\end{abstract}

\noindent
{\bf Keywords}: Parametric model order reduction, parametric dynamical systems, non-intrusive method, minimal rational interpolation, greedy algorithm.
\vspace*{0.5cm}

\noindent
{\bf AMS Subject Classification}: 35B30, 35P15, 41A20, 41A63, 93C35, 93C80.

\section*{Introduction}
The numerical simulation of dynamical systems in frequency domain is of utmost importance in several engineering fields, among which electronic circuit design, acoustics, resonance modeling and control for large structures, and many others. The computational burden of such simulations has kept increasing in the last decades: on one hand, the problem size has been growing because of the need for higher numerical resolution; on the other hand, the necessity to tune design parameters and model uncertain features has lead researchers to tackle parametric models, possibly with a large number of parameters.

The purpose of \emph{model order reduction} (MOR) in general, and of \emph{parametric MOR} (pMOR) in the specific case of dynamical system in the presence of parameters, is to alleviate this computational load. The main strategy to reach this goal relies on building a surrogate model (\emph{reduced order model}, ROM), which mimics accurately the original problem, but which can be solved at a much reduced cost. In the last two decades, the field of pMOR has thrived, leading to the development, analysis, and application of a wide collection of surrogate modeling strategies. In general, we can assign each of these methods to one of two main categories:
\begin{itemize}
\item\textbf{Projection-based pMOR.} The surrogate model is built by restricting the original problem onto a suitable subspace, computed from a set of solutions (most often, snapshots of the system state) of the full problem. This requires access to the operators of the full model, which are not necessarily available in applications, for instance in the case of a black-box solver, or if the system operators never get fully assembled in the solution process. Some subcategories of projective pMOR can be identified depending on whether a global basis (such as POD/Reduced Basis \cite{Baur2011} or multi-parameter multi-moment-matching \cite{Benner2014,Weile1999}) or a collection of local bases (e.g., manifold interpolation of local bases \cite{Amsallem2008} or of reduced system matrices \cite{Amsallem2011,Lohmann2009,Panzer2010}) are employed.
\item\textbf{Non-intrusive pMOR.} The surrogate model is constructed by interpolation or regression of a set of solutions (usually, output samples) of the full problem. As long as the system state is not necessary for the application at hand, it is common to work directly with the system output. In this case as well, the methods can be further split into two subgroups, although the boundary between the two is more vague: some approaches set up the surrogate by solving a unique global interpolation problem \cite{GTalocia2018,Ionita2014,Lefteriu2011}, whereas others build it by first constructing several (rational) models in frequency only, and then combining them over parameter space \cite{Ferranti2011,YueAdapt,YueOverview}.
\end{itemize}

Here, we focus on this last class of techniques: non-intrusive approaches which build a surrogate by interpolating local frequency models, i.e., models constructed at some values of the parameters. A discussion of the pros and cons of this local idea, compared to global approaches, can be found in the excellent survey \cite{Benner2015} (whose focus is, however, mostly on projective pMOR), as well as in \cite{YueOverview}. The latter article is closer in spirit to the present discussion and, as such, we will take it as starting point and main reference for our presentation. Building upon \cite{YueAdapt,YueOverview}, the purpose of this paper is the description of a fully non-intrusive pMOR technique based on parameter interpolation of frequency surrogates. As will be made more clear in the next sections, the objects that are interpolated over parameter space are poles and residues of the Heaviside decomposition of the frequency surrogates. Without going into too much detail here, we summarize briefly the main novelties of our approach, focusing on how it generalizes \cite{YueAdapt}:
\begin{itemize}
\item \emph{Both} frequency space and parameter domain are sampled adaptively, allowing for a better exploration of frequency and parameter domains, as well as for an improved efficiency in the construction of the reduced model. To this aim, we leverage the \emph{minimal rational interpolation} technique as described in \cite{PradoveraGreedy}.
\item Our proposed approach for pole/residue matching has polynomial worst-case complexity despite the combinatorial nature of the task. Notably, our approach can also be applied in the case of unbalanced (i.e., with different numbers of poles and residues) frequency surrogates. This is a critical property in view of our adaptive frequency sampling strategy.
\item We propose a framework for adaptive parameter sampling in a general \emph{high-dimensional} setting by introducing a hierarchical locally-refined sparse grid structure in parameter space.
\end{itemize}
Before proceeding, we deem of importance to remark that local approaches based on interpolation of poles and residues, by their very nature, struggle in dealing with poles and residues that do not depend smoothly on the parameters. In such cases, surrogates that avoid the Heaviside decomposition, e.g., by interpolating the rational frequency response directly rather than its poles and residues \cite{Ferranti2011}, can have more beneficial properties, since they do not rely on the smoothness of poles and residues. Still, we choose to pursue an Heaviside-based approach rather than a rational-based one for two main reasons:
\begin{itemize}
\item By interpolation of poles and residues, one can handle with good flexibility the fairly common case of unbalanced frequency surrogates. More specifically, we will describe how spurious poles may be removed (and missing poles reconstructed) \emph{on the fly} by exploiting the structure of the Heaviside expansion.
\item Despite its intrinsic difficulties in approximating non-smooth poles, an adaptive selection of the parameter sample points, as proposed in this work, can recover a good performance, by refining locally at critical parameter locations (e.g., branch points of the poles).
\end{itemize}

\textbf{Outline of the paper.} We introduce the parametric framework for our approach in Section~\ref{sec:problem}. The ensuing Section~\ref{sec:method} contains our main contribution, in the form of a description of the proposed pMOR technique, with each of the subsections therein examining a different feature of the algorithm. In Section~\ref{sec:bifurcation}, we investigate some of the limitations of our method in the case of crossing and non-smooth poles. We integrate our discussion with two numerical examples in Section~\ref{sec:examples}, showcasing the effectiveness of our method. We conclude with a summary and an outlook for future research in Section~\ref{sec:conclusions}.

\section{Problem framework and notation}\label{sec:problem}
The field of pMOR is closely entwined with the study of parametric dynamical systems in the frequency domain. More precisely, one considers the problem:
\begin{equation}\label{eq:dynamical}
\text{given }(z,\pp)\in Z\times\PP\text{, find }Y=Y(z,\pp)\text{ such that }\begin{cases}
\left(zE_\pp-A_\pp\right)X(z,\pp)=B_\pp,\\
Y(z,\pp)=C_\pp X(z,\pp),
\end{cases}
\end{equation}
where
\begin{equation*}
A_\pp,E_\pp\in\xC^{n_S\times n_S},\quad X(z,\pp),B_\pp\in\xC^{n_S\times n_I},\quad C_\pp\in\xC^{n_O\times n_S},\quad\text{and}\quad Y(z,\pp)\in\xC^{n_O\times n_I}.
\end{equation*}
We call $X$ and $Y$ the \emph{state} and the \emph{output} of the system, respectively. The sets $Z\subset\xC$ and $\PP\subset\xC^d$ ($\xR^d$ in most applications) are frequency and parameter domains, respectively. The subscript $\pp$ of the system matrices denotes their eventual dependence on the parameters $\pp$. We remark that we do not exclude the case of the state being the output of interest: such case can be obtained quite trivially by setting $n_O=n_S$ and $C_\pp$ equal to the identity matrix.

The behavior of state $X$ and output $Y$ with respect to $z$ only (for fixed $\pp$, the non-parametric case) is well understood, and is backed up by an extensive literature in system theory \cite{Antoulas2017}. Among the many properties of these systems, one is crucial to our discussion, namely the Heaviside decomposition: under some quite broad assumptions on the spectral properties of the pencil $(A_\pp,E_\pp)$, we can write
\begin{equation}\label{eq:output}
Y(z,\pp)=\sum_{j}(z-\lambda_\pp^{(j)})^{-m_\pp^{(j)}}Y_\pp^{(j)}\quad\text{for }z\in Z\setminus\{\lambda_\pp^{(j)}\}_{j},
\end{equation}
with $\lambda_\pp^{(j)}\in\xC\cup\{\infty\}$, $m_\pp^{(j)}\in\xN$, and $Y_\pp^{(j)}\subset\xC^{n_O\times n_I}$, for all $j$.

We remark that our discussion relies only on the decomposition \eqref{eq:output} of the output, together with the following assumptions:
\begin{enumerate}[(a)]
\item $Y:Z\times\PP\to V$, with $V$ a normed vector space (for \eqref{eq:dynamical}, we may set $V=\xC^{n_O\times n_I}$ endowed with the Frobenius norm).
\item In \eqref{eq:output}, $\lambda^{(j)}$ and $Y^{(j)}$ depend continuously on $\pp$ for all $j$, cf. Section~\ref{sec:bifurcation}.
\item\label{item:multiplicity} In \eqref{eq:output}, $m^{(j)}$ is independent of $\pp$ for all $j$.
\end{enumerate}
As such, our proposed strategy extends further than linear parametric dynamical systems \eqref{eq:dynamical}. Some examples of practical interest are parametric scattering problems in frequency domain \cite{Hiptmair2018} and parametric nonlinear eigenproblems \cite{Betcke2013}. Still, for simplicity of exposition, we will restrict our discussion to the finite-dimensional linear parametric dynamical system \eqref{eq:dynamical}.

In order to simplify the presentation, from here onward we will replace \eqref{item:multiplicity} with the stronger assumption that $m^{(j)}$ be equal to 1. This excludes the possibility of degenerate eigenvalues. We remark that, in most practical applications, each multiplicity $m^{(j)}$ is indeed identically equal to 1. We postpone a discussion on this till Sections~\ref{sec:bifurcation} and \ref{sec:conclusions}.

\section{The double-greedy pMOR strategy}\label{sec:method}
We start this section by detailing our pMOR technique in its formulation without $\pp$-adaptivity. The general structure of the pole-matching-based pMOR algorithm is summarized in Algorithms \ref{algo:build} and \ref{algo:evaluate}, and the different building blocks are discussed more extensively in the following subsections.

An online-offline decomposition of the algorithm is performed, in the usual MOR fashion \cite{Quarteroni2015}: first, the surrogate model is built in an expensive training phase, which requires solving the original problem at several values of frequency and parameters; the reduced model is then stored and can be evaluated with a (hopefully) much reduced computational cost at arbitrary frequency and parameter values.

\begin{algorithm}[t]
    \caption{Pole-matching pMOR -- Offline phase}
    \label{algo:build}
    \begin{algorithmic}[1]
    		\Function{pROM\_Train}{$\pp_1,\ldots,\pp_S$}
    		\For{$k=1,\ldots,S$}
	    		\State $\textup{ROM}_k$ $\gets$ \Call{BuildFrequencyROM}{$\pp_k$} \Comment{Section~\ref{sec:zadapt}}
    		\EndFor
    		\State $J\gets\big\{1\big\}$ \Comment{Initialize root of search tree}
    		\For{$\textrm{iter}=2,\ldots,S$} \Comment{Matching loop}
    			\State $(k,\ell)\gets\argmin_{k'\in J,\,\ell'\notin J}\norm{\pp_{k'}-\pp_{\ell'}}$ \Comment{Breadth-first search}
    			\State permutation $\gets$ \Call{ModalMatching}{$\textup{ROM}_k$, $\textup{ROM}_\ell$}\Comment{Section~\ref{sec:matching}}
	    		\State $\textup{ROM}_\ell\gets$\Call{ApplyPermutation}{$\textup{ROM}_\ell$, permutation}
	    		\State $J\gets J\cup\{\ell\}$ \Comment{Add to explored set}
    		\EndFor
    		\State $\ppsi\gets$\Call{BuildLocalWeightFunctions}{$\pp_1,\ldots,\pp_S$} \Comment{Section~\ref{sec:interp}}
    		\State \textbf{return} $\{\textup{ROM}_k\}_{k=1}^S$, $\ppsi$
    		\EndFunction
    \end{algorithmic}
\end{algorithm}

The offline phase is shown in Algorithm~\ref{algo:build}. The input of the procedure is a set of parameter values $\PP_{train}=\{\pp_1,\ldots,\pp_S\}\subset\PP$, where reduced models \emph{in frequency only} are built: more precisely, for each $k=1,\ldots,S$, we compute the approximation
\begin{equation}\label{eq:ROMs}
Y(z,\pp_k)\approx\textup{ROM}_k(z)=\sum_{j=1}^{\overline{n}_k}\frac{\overline{Y}_k^{(j)}}{z-\overline{\lambda}_k^{(j)}}.
\end{equation}
More details on this step are given in Section~\ref{sec:zadapt}. The addends of the sum in \eqref{eq:ROMs} are then sorted in such a way that poles $\overline{\lambda}_k^{(j)}$ and residues $\overline{Y}_k^{(j)}$ with the same $j$ but different $k$ ``correspond to each other''. We explain what we mean by this, and how we achieve it, in Section~\ref{sec:matching}. In the same section, we also discuss how we deal with the situation where two surrogate models have different amounts of poles, i.e., $\overline{n}_k\neq\overline{n}_\ell$. For the remainder of the present overview, for simplicity we assume $\overline{n}_1=\ldots=\overline{n}_S=:\overline{n}$.

\begin{algorithm}[t]
    \caption{Pole-matching pMOR -- Online phase}
    \label{algo:evaluate}
    \begin{algorithmic}[1]
    		\Function{pROM}{$z$, $\pp$; $\{\textup{ROM}_k\}_{k=1}^S$, $\ppsi$}
    		\State $\widetilde{Y}\gets 0$,\ $\mathbf{w}\gets\ppsi(\pp)$
    		\For{$j=1,2,\ldots,\overline{n}$} \Comment{Loop over surrogate poles}
	    		\State $\widetilde{\lambda}^{(j)}\gets 0$, $\widetilde{Y}^{(j)}\gets 0$
	   		\For{$k=1,2,\ldots,S$} \Comment{Loop over frequency models}
		    		\State Extract $\overline{\lambda}_k^{(j)}$ and $\overline{Y}_k^{(j)}$ from $\textup{ROM}_k$
	    			\State $\widetilde{\lambda}^{(j)}\gets \widetilde{\lambda}^{(j)}+w_k\overline{\lambda}_k^{(j)}$,\ \ $\widetilde{Y}^{(j)}\gets \widetilde{Y}^{(j)}+w_k\overline{Y}_k^{(j)}$ \Comment{Add local contribution}
	    		\EndFor
    			\State $\widetilde{Y}\gets\widetilde{Y}+\widetilde{Y}^{(j)}/\big(z-\widetilde{\lambda}^{(j)}\big)$ \Comment{Add Heaviside term}
    		\EndFor
    		\State \textbf{return} $\widetilde{Y}$
    		\EndFunction
    \end{algorithmic}
\end{algorithm}

At this point, it only remains to prescribe a rule to define the reduced model at a new parameter value $\pp\in\PP\setminus\PP_{train}$. To this aim, we define $S$ weight functions $\ppsi=(\psi_1,\ldots,\psi_S):\PP\to\xC^S$, which we employ to interpolate poles and residues over $\pp$-space:
\begin{equation}\label{eq:interpolate}
\widetilde{\lambda}^{(j)}(\pp)=\sum_{k=1}^S\psi_k(\pp)\overline{\lambda}_k^{(j)}\quad\text{and}\quad\widetilde{Y}^{(j)}(\pp)=\sum_{k=1}^S\psi_k(\pp)\overline{Y}_k^{(j)}.
\end{equation}
The resulting global model can then be evaluated through Algorithm~\ref{algo:evaluate}, as
\begin{equation}\label{eq:ROMglobal}
Y(z,\pp)\approx\textup{pROM}(z,\pp)=\sum_{j=1}^{\overline{n}}\frac{\widetilde{Y}^{(j)}(\pp)}{z-\widetilde{\lambda}^{(j)}(\pp)}=\sum_{j=1}^{\overline{n}}\frac{\sum_{k=1}^S\psi_k(\pp)\overline{Y}_k^{(j)}}{z-\sum_{k=1}^S\psi_k(\pp)\overline{\lambda}_k^{(j)}}.
\end{equation}
Note that, in order for the frequency surrogates to be interpolated exactly, the weights should satisfy the conditions
\begin{equation}\label{eq:weights}
\psi_\ell(\pp_k)=\begin{cases}
1\quad\text{if }k=\ell,\\
0\quad\text{if }k\neq\ell,
\end{cases}
\end{equation}
cf. Section~\ref{sec:interp}.

\begin{rmrk}\label{rem:polynomial}
Depending on the choice of MOR method in frequency, the Heaviside formulation \eqref{eq:ROMs} might also include constant, polynomial, or, more generally, smooth terms: for instance,
\begin{equation}\label{eq:ROMspoly}
\textup{ROM}_k(z)=\sum_{j=1}^{\overline{n}_k}\frac{\overline{Y}_k^{(j)}}{z-\overline{\lambda}_k^{(j)}}+\sum_{\ell=0}^{\overline{m}}\overline{W}_k^{(\ell)}z^\ell.
\end{equation}
If this is the case, we do not need to modify our algorithm, because the additional terms do not require to be matched nor permuted. The resulting global surrogate (after pole matching) is simply
\begin{equation*}
\textup{pROM}(z,\pp)=\sum_{j=1}^{\overline{n}}\frac{\sum_{k=1}^S\psi_k(\pp)\overline{Y}_k^{(j)}}{z-\sum_{k=1}^S\psi_k(\pp)\overline{\lambda}_k^{(j)}}+\sum_{\ell=0}^{\overline{m}}\sum_{k=1}^S\psi_k(\pp)\overline{W}_k^{(\ell)}z^\ell.
\end{equation*}
\end{rmrk}

\subsection{Frequency adaptivity via Minimal Rational Interpolation}\label{sec:zadapt}
Here we provide some details on the function \textsc{BuildFrequencyROM} in Algorithm~\ref{algo:build}, which encodes the construction of a ROM for a non-parametric dynamical system with respect to a single parameter, namely the frequency $z$. A number of surrogate modeling strategies for such problems have been proposed in the MOR literature: the most famous are, among projection-based methods, the Proper Orthogonal Decomposition/Reduced Basis \cite{Benner2015} and the Krylov/Moment Matching \cite{Benner2014} methods, and, among the non-intrusive techniques, the Loewner Framework \cite{Ionita2014} and the Vector Fitting algorithm \cite{Drmac2015,Grivet2016}. Any of them could be used to supply the surrogate modeling that we require. In fact, different frequency surrogates could even be obtained by different MOR approaches, as long as each reduced model allows for a Heaviside expansion \eqref{eq:ROMs}.

\begin{algorithm}[t]
    \caption{Greedy Minimal Rational Interpolation}
    \label{algo:mri}
    \begin{algorithmic}[1]
    		\Function{BuildFrequencyROM}{$\pp$}
    		\State Initialize $Z_{train}\subset Z$ (very coarse), $Z_{test}\subset Z$ (fine), $\text{tol}>0$
    		\State snapshots $\gets\{\,\}$
    		\For{$z\in Z_{train}$} \Comment{Model initialization}
    			\State Evaluate $Y(z,\pp)$ \Comment{Full model solve}
	    		\State snapshots $\gets$ snapshots $\cup\ \{Y(z,\pp)\}$
    		\EndFor
    		\Repeat \Comment{Greedy loop}
	    		\State ROM$\gets$ \Call{BuildMRI}{$Z_{train}$, snapshots} \Comment{SVD of a matrix of size $\#Z_{train}$\cite{Pradovera2020}}
	    		\State $z^\star\gets$ \Call{GetNextSamplePoint}{ROM, $Z_{test}$} \Comment{Maximize \eqref{eq:aposteriori} over $Z_{test}$}
    			\State Evaluate $Y(z^\star,\pp)$ \Comment{Full model solve}
	    		\State snapshots $\gets$ snapshots $\cup\ \{Y(z^\star,\pp)\}$
	    		\State Move $z^\star$ from $Z_{test}$ to $Z_{train}$
    		\Until{$\norm{Y(z^\star,\pp)-\textup{ROM}(z^\star)}_V\leq$ tol $\norm{Y(z^\star,\pp)}_V$}
    		\State \emph{Optional}:\ ROM$\gets$ \Call{BuildMRI}{$Z_{train}$, snapshots} \Comment{Remark~\ref{rem:nowaste}}
    		\State \emph{Optional}:\ ROM$\gets$ \Call{CleanUp}{ROM} \Comment{Remark~\ref{rem:cleanup}}
    		\State \textbf{return} ROM
    		\EndFunction
    \end{algorithmic}
\end{algorithm}

In this work, we consider the Minimal Rational Interpolation (MRI) method proposed in \cite{Pradovera2020}, generalizing \cite{Bonizzoni}. MRI is non-intrusive, i.e., it does not require access to the matrices appearing in the dynamical system \eqref{eq:dynamical}, allowing for a wider applicability of the method. In particular, MRI can be used to efficiently build surrogates for vector-valued quantities, e.g., high-(or even $\infty$-)dimensional states $X$ of dynamical systems \eqref{eq:dynamical}. (In fact, the effectiveness of MRI may even improve if the ambient space, where snapshots of the approximation target are located, is large.) At the same time, one can apply it in a greedy fashion \cite{PradoveraGreedy}, so that the number and location of frequency samples are selected adaptively, in such a way that a prescribed accuracy is attained over the whole frequency domain $Z$. The specific greedy strategy that we choose has a proper theoretical motivation only if the parametric problem depends on $z$ in a simple way, e.g., linearly, as in \eqref{eq:dynamical}. If this is not the case, or if the dependence on frequency is not known, one may want to consider possible alternatives \cite{PradoveraGreedy}. More details on the method are provided in \cite{Pradovera2020}. Here we only give a short overview, which is summarized in Algorithm~\ref{algo:mri}.

Let a parameter $\pp\in\PP$ be fixed. Given $S$ distinct frequency sample points $Z_{train}\subset Z$ and corresponding snapshots $\{Y(z,\pp)\}_{z\in Z_{train}}$, the MRI procedure first builds a surrogate denominator $Q$ of degree $\leq S-1$ by solving a minimization problem involving the snapshots: in practice, this problem can be solved at $\mathcal{O}(S^3)$ computational cost by SVD. Now the surrogate poles $\{\overline{\lambda}^{(j)}\}_{j=1}^{\overline{n}}$ ($\overline{n}\leq S-1$ is the degree of $Q$) appearing in \eqref{eq:ROMspoly} can be extracted from $Q$ through any root-finding algorithm\footnote{We are assuming all poles to be simple, cf. Section~\ref{sec:ex1b}.}. Then we set $\overline{m}=S-\overline{n}-1$, and we compute the $S$ terms $\{\overline{Y}^{(j)}\}_{j=1}^{\overline{n}}\cup\{\overline{W}^{(\ell)}\}_{\ell=0}^{\overline{m}}$ in the ROM \eqref{eq:ROMspoly} by interpolation:
\begin{equation}\label{eq:MRIinterp}
\sum_{j=1}^{\overline{n}}\frac{\overline{Y}^{(j)}}{z-\overline{\lambda}^{(j)}}+\sum_{\ell=0}^{\overline{m}}\overline{W}^{(\ell)}z^\ell=Y(z,\pp)\qquad\forall z\in Z_{train}.
\end{equation}
The procedure just described is encoded by the function \textsc{BuildMRI} in Algorithm~\ref{algo:mri}.

The adaptive selection of frequency sample points is carried out via the typical greedy-MOR loop: at each iteration a new sample gets added, at a position selected in the test set $Z_{test}$, based on the current reduced model. More precisely, the next sample $z^\star$ is selected as the maximizer over $Z_{test}$ of some \emph{a posteriori} indicator, and the algorithm terminates when such indicator is smaller than a prescribed tolerance. However, in our framework, standard residual estimators (e.g., the classic Reduced Basis one \cite{Quarteroni2015}) cannot be employed in a non-intrusive manner. Instead, we apply the ``look-ahead'' idea introduced in \cite{PradoveraGreedy}, to which we refer for more details:
\begin{itemize}
\item Exploiting only the current reduced model, we select the next sample point $z^\star$, as the maximizer over $Z_{test}$ of
\begin{equation}\label{eq:aposteriori}
z\mapsto\prod_{z'\in Z_{train}}\abs{z-z'}\ \bigg/\ \prod_{j=1}^{\overline{n}}\abs{z-\overline{\lambda}^{(j)}}.
\end{equation}
\item Since an unknown scaling constant is involved, \eqref{eq:aposteriori} does not tell us whether the prescribed tolerance is satisfied; to obtain this information, we perform an expensive solve of the full model at $z^\star$, and evaluate the approximation error explicitly.
\end{itemize}

\begin{rmrk}\label{rem:nowaste}
As is evident from Algorithm~\ref{algo:mri}, the extra expensive solve of the full system does not go to waste, as it is exactly the snapshot which is necessary at the next iteration. Only at the last iteration, when the tolerance is finally satisfied, the additional solution is not included in the surrogate. In practice, to avoid wasting the ``free'' final snapshot, we can actually add even this last sample to the reduced model by running \textsc{BuildMRI} an additional time after the greedy loop is completed.
\end{rmrk}

\begin{rmrk}\label{rem:cleanup}
As discussed in more detail in Section~\ref{sec:matching}, for a reliable matching of poles and residues of different frequency \textup{ROMs}, it is crucial to remove any unwanted spurious (sometimes referred to as ``parasitic'' or ``non-physical'') pole-residue pairs. As such, it is usually worth the effort to add a post-processing ``clean-up'' step, to remove unwanted pole-residue pairs. The development of reliable strategies to identify spurious effects remains an open problem in rational approximation. However, two quite simple and inexpensive criteria that can be applied are:
\begin{itemize}
\item Remove poles which are too far away from $Z$.
\item Remove poles whose residues have small magnitude (possible Froissart doublets \cite{Beckermann2018}).
\end{itemize}
In addition, it might be possible to remove poles which are provably spurious based on the properties of the problem, e.g., non-real poles for self-adjoint problems, or unstable poles for stable systems.

Any pole removed in this context should be automatically balanced by increasing the number of smooth terms in the ROM, namely $\overline{m}$ in Remark~\ref{rem:polynomial}. We refer to \cite{Pradovera2020} for a more detailed discussion on this procedure, which, employing the notation therein, could be summarized as ``choosing $N<S-1$''. Thanks to this strategy, if the eliminated poles are indeed spurious, their removal \emph{does not} impact negatively the accuracy of the frequency surrogate.
\end{rmrk}

At this point, it is important to note that, in general, we cannot guarantee all the frequency ROMs to have the same number of pole-residue pairs $\overline{n}$, even after a perfect removal of all the spurious effects. The main reason for this is that the full model \eqref{eq:dynamical} may have a different number of poles in $Z$ for different values of $\pp$. A numerical example showcasing poles entering and leaving the frequency domain can be found in Section~\ref{sec:ex2}. This can quickly become problematic, since the greedy MRI algorithm is not guaranteed to identify well poles outside $Z$ (actually, even to capture them at all). One could try to counteract this issue by building frequency models on an enlarged domain $Z'\supset Z$, at the cost of a higher offline time. Still, the same issue of migrating poles might simply arise again at a larger scale.

\subsection{Matching frequency models}\label{sec:matching}
The offline matching loop in Algorithm~\ref{algo:build} takes care of permuting the Heaviside addends of the frequency ROMs, see \eqref{eq:ROMs}. The objective of this step is to identify with good accuracy how the location of each pole-residue pair evolves as $\pp$ changes: this means that, for each $j$, we wish $\overline{\lambda}_1^{(j)},\overline{\lambda}_2^{(j)},\ldots,\overline{\lambda}_S^{(j)}$ (resp. $\overline{Y}_1^{(j)},\overline{Y}_2^{(j)},\ldots,\overline{Y}_S^{(j)}$) to be approximations of the ``same'' pole (resp. residue) of the parametric problem for $\pp=\pp_1,\pp_2,\ldots,\pp_S$. By ``same'' pole (resp. residue) we mean some continuous curve $\pp\mapsto\lambda_\pp^{(j')}$ (resp. $\pp\mapsto Y_\pp^{(j')}$) for some fixed $j'$, with $\lambda_\pp^{(j')}$ and $Y_\pp^{(j')}$ defined in Section~\ref{sec:problem}. Note that this notion is not unique whenever poles intersect. This case is discussed in Section~\ref{sec:bifurcation}.

The frequency surrogate models are explored in a breadth-first fashion, matching one new ROM at a time to its closest (with respect to $\pp$) neighbor. We exclude a global approach, matching all frequency models at the same time, due to its computational unfeasibility: $S$-partite matching is notoriously an NP-hard problem for $S\geq 3$ \cite{Karp1972}.

Let us characterize one of this bipartite matching problems, namely the one between $\textup{ROM}_k$ and $\textup{ROM}_\ell$. For the remainder of this section, we assume that the two models have the same number of poles, i.e., $\overline{n}_k=\overline{n}_\ell=:\overline{n}$. The more general and, due to our choice of frequency surrogate modeling, (potentially) more common case $\overline{n}_k\neq\overline{n}_\ell$ is discussed in Section~\ref{sec:unbalanced}. We seek $\ssigma=(\sigma_1,\ldots,\sigma_{\overline{n}})$ which attains
\begin{equation}\label{eq:matchfunc}
\min_{\ssigma\in(1,2,\ldots,\overline{n})!}\left(\sum_{j=1}^{\overline{n}}\abs{\overline{\lambda}_k^{(j)}-\overline{\lambda}_\ell^{(\sigma_j)}}+w\sum_{j=1}^{\overline{n}}\norm{\overline{Y}_k^{(j)}-\overline{Y}_\ell^{(\sigma_j)}}_V\right),
\end{equation}
where $(1,2,\ldots,\overline{n})!$ denotes the set of permutations of the tuple $(1,2,\ldots,\overline{n})$. We will refer to the minimal value of \eqref{eq:matchfunc} as the \emph{Heaviside distance}\footnote{When smooth terms are added, see \eqref{eq:ROMspoly}, or in the unbalanced case $\overline{n}_k\neq\overline{n}_\ell$, see \eqref{eq:matchfuncsup}, this is not an actual distance because it is only semi-positive definite.} between models $\textup{ROM}_k$ and $\textup{ROM}_\ell$. In \eqref{eq:matchfunc}, $w\in[0,\infty]$ denotes a weight, representing the relative importance given to poles and residues in the matching. We note that a cost functional almost identical to \eqref{eq:matchfunc} was first introduced in \cite{YueAdapt} to tackle the pole-residue matching.

An alternative formulation of the same problem is: given the cost matrix $D\in\xR^{\overline{n}\times\overline{n}}$, with entries
\begin{equation}\label{eq:matchcost}
D_{j,j'}=\abs{\overline{\lambda}_k^{(j)}-\overline{\lambda}_\ell^{(j')}}+w\norm{\overline{Y}_k^{(j)}-\overline{Y}_\ell^{(j')}}_V,
\end{equation}
we wish to extract exactly one value per row and one value per column, so that the sum of the selected entries is minimized. Despite the combinatorial nature of this optimization problem, there exist polynomial-time algorithms to solve it (e.g., using ideas from maximum flow problems \cite{Crouse2016}). Here, we solve this problem via the \texttt{linear\_sum\_assignment} function in the \texttt{scipy.optimize} module \cite{SciPy}, which requires as only input the cost matrix $D$. As shown in \cite[Section II.C]{Crouse2016}, this implementation has $\mathcal{O}(\overline{n}^3)$ worst-case complexity, but we note that only $\mathcal{O}(\overline{n}^2)$ operations are necessary if just $\mathcal{O}(1)$ pole swaps are needed, i.e., if the optimal $\ssigma$ is ``close'' to $(1,2,\ldots,\overline{n})$.

\begin{rmrk}\label{rem:order}
When $\overline{n}_k=\overline{n}_\ell$, the matching optimization problem is symmetric. By this, we mean that, once the optimal permutation has been found, we can either apply it to the ``new'' model $\textup{ROM}_\ell$, or apply its inverse to the ``old'' model $\textup{ROM}_k$. In practice, in the scope of the breadth-first search in Algorithm~\ref{algo:build}, it is convenient to always rearrange the poles-residues of the new model $\ell$, since applying the permutation to model $k$ would require permuting \emph{all} the already explored models $\{\textup{ROM}_{k'}\}_{k'\in J}$.
\end{rmrk}

\subsubsection{Unbalanced matching}\label{sec:unbalanced}
Suppose $\overline{n}_k<\overline{n}_\ell$ (the converse case can be treated analogously by considering a ``transposed'' matching problem). We can cast the matching optimization problem in rectangular form: find $\ssigma=(\sigma_1,\ldots,\sigma_{\overline{n}_\ell})$ which attains
\begin{equation}\label{eq:matchfuncsup}
\min_{\ssigma\in(1,2,\ldots,\overline{n}_\ell)!}\left(\sum_{j=1}^{\overline{n}_k}\abs{\overline{\lambda}_k^{(j)}-\overline{\lambda}_\ell^{(\sigma_j)}}+w\sum_{j=1}^{\overline{n}_k}\norm{\overline{Y}_k^{(j)}-\overline{Y}_\ell^{(\sigma_j)}}_V\right).
\end{equation}
(The minimization of \eqref{eq:matchfuncsup} does not require more computational effort than a balanced problem \eqref{eq:matchfunc} of size $\overline{n}_\ell$ \cite{Crouse2016}.) Then $(\sigma_1,\ldots,\sigma_{\overline{n}_k})$ gives the desired permutation, while the indices $\sigma_{\overline{n}_k+1},\ldots,\sigma_{\overline{n}_\ell}$ remain unassigned. It remains to choose how to deal with the $\overline{n}_\ell-\overline{n}_k$ extra pole-residue pairs with indices $\sigma_{\overline{n}_k+1},\ldots,\sigma_{\overline{n}_\ell}$: this choice depends on whether we think that they are ($\circ$) missing from $\textup{ROM}_k$ or ($\bullet$) spurious in $\textup{ROM}_\ell$.

The solution for case \textup{($\bullet$)} is quite straightforward: we simply remove the $\overline{n}_\ell-\overline{n}_k$ erroneous poles and residues from $\textup{ROM}_\ell$, as well as from all the surrogates matched to it. Note, however, that the $\overline{n}_k$ remaining residues of $\textup{ROM}_\ell$ should be updated after the removal of the spurious poles, in order to guarantee good approximation properties at $\pp_\ell$. To this aim, it is quite cheap to recompute the residues of $\textup{ROM}_\ell$ by solving an interpolation problem, cf. \eqref{eq:MRIinterp}, exploiting the already available snapshots at $\pp_\ell$. In particular, the approach described in Remark~\ref{rem:cleanup}, i.e., the addition of smooth terms to compensate for the removed poles/residues, can be applied here.

Instead, in case \textup{($\circ$)}, the problem is much harder: we wish to reconstruct the poles and residues unaccounted for at $\pp_k$ from information at $\pp_\ell$. One naive way to achieve this involves ``copying'' the missing poles and residues from $\textup{ROM}_\ell$ to $\textup{ROM}_k$, i.e., appending the synthetic terms
\begin{equation}\label{eq:matchunbalanced}
\frac{\overline{Y}_\ell^{(\sigma_{\overline{n}_k+1})}}{z-\overline{\lambda}_\ell^{(\sigma_{\overline{n}_k+1})}}+\ldots+\frac{\overline{Y}_\ell^{(\sigma_{\overline{n}_\ell})}}{z-\overline{\lambda}_\ell^{(\sigma_{\overline{n}_\ell})}}
\end{equation}
to the Heaviside expansion of $\textup{ROM}_k$. This results in an adjusted surrogate, with a balanced pole matching. Note that, similarly to pole removal, the $\overline{n}_k$ original residues of $\textup{ROM}_k$ should be updated to account for the added synthetic terms. We remark that more refined (but also potentially less stable) approaches may be based on other forms of extrapolation of the additional Heaviside terms, not only from $\textup{ROM}_\ell$, but also from all the other frequency surrogates that contain the poles with indices $\sigma_{\overline{n}_k+1},\ldots,\sigma_{\overline{n}_\ell}$, cf. Remark~\ref{rem:reconstruct} and Section~\ref{sec:interp}.

\begin{rmrk}\label{rem:retrain}
In case \textup{($\circ$)}, one could impose a balanced matching by retraining the poorer model $\textup{ROM}_k$, forcing more iterations of greedy MRI until the surrogate has exactly $\overline{n}_\ell$ poles. However, this approach has two clear issues, which might make it disadvantageous in practice:
\begin{itemize}
\item If many surrogates need retraining, one may incur in substantial additional computational cost. This is the case especially if $\pp$-adaptivity is employed, see Section~\ref{sec:padapt}.
\item Instead of identifying correctly the missing poles, MRI could introduce spurious effects, thus interfering with the matching procedure, rather than helping it. This is likely to happen if the missing poles are located outside the frequency domain $Z$.
\end{itemize}
\end{rmrk}

At this point, we deem important to give a \emph{caveat}: neither of the two approaches above (pole reconstruction and removal) is, on its own, able to solve adequately all practical situations, as we showcase in a practical example in Section~\ref{sec:ex2}. Instead, a hybrid version, with some poles getting removed and some reconstructed, has the potential to perform better than either approach. Here, we choose to pursue this third approach, whose effectiveness obviously hinges on how well we can differentiate between ``good'' and ``bad'' poles.

\begin{algorithm}[t]
    \caption{General pole-matching}
    \label{algo:unbalanced}
    \begin{algorithmic}[1]
    		\State $J\gets\big\{1\big\}$ \Comment{Initialize root of search tree, Remark~\ref{rem:root}}
    		\For{$\textrm{iter}=2,\ldots,S$} \Comment{Matching loop}
    			\State $(k,\ell)\gets\argmin_{k'\in J,\,\ell'\notin J}\norm{\pp_{k'}-\pp_{\ell'}}$ \Comment{Breadth-first search}
    			\State $\ssigma\gets$ \Call{ModalMatching}{$\textup{ROM}_k$, $\textup{ROM}_\ell$} \Comment{\eqref{eq:matchfunc} or \eqref{eq:matchfuncsup}}
    			\State $\textup{ROM}_\ell\gets$\Call{ApplyPermutation}{$\textup{ROM}_\ell$, $\ssigma$} \Comment{Possibly using only part of $\ssigma$}
    			\If{$\overline{n}_k>\overline{n}_\ell$}
	    			\State $\textup{ROM}_\ell\gets\textup{ROM}_\ell+\sum_{j=\overline{n}_\ell+1}^{\overline{n}_k}\overline{Y}_k^{(j)}/\big(\cdot-\overline{\lambda}_k^{(j)}\big)$\Comment{Extrapolate at index $\ell$}
		    		\State Flag the $\overline{n}_k-\overline{n}_\ell$ poles added to $\textup{ROM}_\ell$ as ``synthetic''
	    			\State $\overline{n}_\ell\gets\overline{n}_k$
    			\ElsIf{$\overline{n}_k<\overline{n}_\ell$}
		    		\For{$i\in J$}\Comment{Extrapolate at indices $J$}
		    			\State $\textup{ROM}_i\gets\textup{ROM}_i+\sum_{j=\overline{n}_i+1}^{\overline{n}_\ell}\overline{Y}_\ell^{(j)}/\big(\cdot-\overline{\lambda}_\ell^{(j)}\big)$			    	
		    			\State Flag the $\overline{n}_\ell-\overline{n}_i$ poles added to $\textup{ROM}_i$ as ``synthetic''
		    			\State $\overline{n}_i\gets\overline{n}_\ell$
	    			\EndFor
    			\EndIf
	    		\State $J\gets J\cup\{\ell\}$ \Comment{Add to explored set}
    		\EndFor
    		\For{$j=1,\ldots,\max_{j\in J}\overline{n}_j$}
    			\State synthetic\_count $\gets\#\{k=1,2,\ldots,S$ such that $\overline{\lambda}_k^{(j)}$ is synthetic$\}$
    			\If{synthetic\_count $>S(1-$ tol\textsubscript{synth}$)$}
    				\State Remove $\overline{\lambda}^{(j)}$ from all ROMs \Comment{Remove poles which are synthetic too many times}
    			\EndIf
    		\EndFor
    		\State \emph{Optional}: apply higher-order reconstruction to all remaining synthetic poles \Comment{Remark~\ref{rem:reconstruct}}
    \end{algorithmic}
\end{algorithm}

We summarize our proposed strategy in Algorithm~\ref{algo:unbalanced}, which replaces the simple matching loop in Algorithm~\ref{algo:build}. The main idea is the following: whenever it becomes necessary to match unbalanced models, the naive reconstruction \eqref{eq:matchunbalanced} is used to augment the less rich surrogate(s); however, the added synthetic poles are flagged as unreliable. At the end of the matching loop, all the models contain the same number of poles $\overline{n}=\max_{k=1,\ldots,S}\overline{n}_k$. At this point, if a pole with a certain index is too often unreliable, it is removed from the pROM. Here, ``too often'' is determined based on a given tolerance tol\textsubscript{synth} between $0$ and $1$: the extreme values $0$ and $1$ correspond to cases \textup{($\circ$)} and \textup{($\bullet$)}, respectively.

\begin{rmrk}\label{rem:root}
\begin{figure}[t!]
\centering
\includegraphics{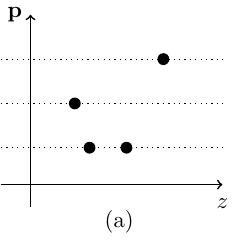}\hspace{.5cm}%
\includegraphics{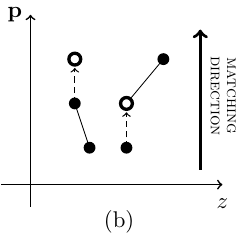}\hspace{.5cm}%
\includegraphics{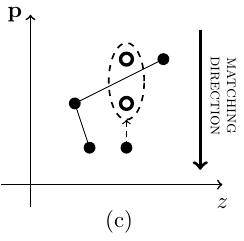}
\caption{Example of history-dependent matching. Poles of the surrogates are denoted by full dots (a). Results by matching bottom-up (b) and top-down (c). Dashed arrows and empty circles are used to denote pole reconstruction steps and synthetic poles, respectively.}
\label{fig:matchingorder}
\end{figure}

The reconstruction of missing poles introduces an asymmetry in the matching procedure, so that the order in which the models are matched matters: for instance, see the situation depicted in Figure~\ref{fig:matchingorder}. As such, it seems important to choose well the root in the breadth-first exploration of $\PP_{train}$. However, choosing optimally this root remains an open problem. From a computational point of view, it makes sense to choose as root the surrogate with the largest number of poles, so that we never have to retrace our steps to add synthetic poles.
\end{rmrk}

\begin{rmrk}\label{rem:reconstruct}
The constant pole reconstruction \eqref{eq:matchunbalanced} is quite blunt, especially when the parameter resolution is low. If the poles depend smoothly on $\pp$, it is preferable to employ a reconstruction with a larger stencil, for instance global (least squares) polynomial extrapolation, using information from all the surrogates which contain the missing pole, see Section~\ref{sec:interp}. However, this is not always viable during the matching loop, since we explore $\PP_{train}$ breadth-first: for instance, if $\overline{n}_k<\overline{n}_\ell$, we are forced to reconstruct poles from the single new model $\textup{ROM}_\ell$ (a similar problem may arise in the case $\overline{n}_k>\overline{n}_\ell$ if $J$ is too small). Still, as shown in the last line of Algorithm~\ref{algo:unbalanced}, it remains feasible to apply a higher-order pole reconstruction after the matching loop is complete.
\end{rmrk}

\subsection{Global vs local interpolation}\label{sec:interp}
In this section we describe how one can employ the local pole-residue information \eqref{eq:ROMs} at $\PP_{train}$ to obtain an approximate Heaviside expansion at a new point $\pp\in\PP\setminus\PP_{train}$. As shown in \eqref{eq:interpolate}, we rely on an interpolation strategy encoded by the weight functions $\ppsi:\PP\to\xC^S$, so that we simply need to prescribe how such weights are constructed. It is important to note that this task, in its natural formulation, is independent of the frequency surrogates, and depends only on the location of the sample parameter points $\PP_{train}$, cf. Remark~\ref{rem:passivity}.

Setting \eqref{eq:weights} as target, we can cast the problem of finding $\ppsi$ as $S$ independent $\PP\to\xC$ interpolation problems (one for each $\psi_\ell$). To address this problem, an extensive amount of techniques and results are available in the literature \cite{Chkifa2014,Phillips2003,Wendland2004}. Here, we consider 3 options:
\begin{itemize}
\item \textbf{Global.} We can seek weights $\ppsi$ within some function space with global regularity, e.g., polynomials, or radial basis functions (using a smooth kernel to achieve interpolation). This approach can potentially achieve high accuracy, but relies on some level of smoothness of poles and residues with respect to $\pp$. It is important to note that, in this framework, the interpolation condition \eqref{eq:weights} could be weakened and enforced only in a least-squares sense: the resulting surrogate would then require less memory for the storage of approximate poles and residues, cf. the ``regression'' step in \cite{YueAdapt}.
\item \textbf{Local structured.} To satisfy \eqref{eq:weights}, we can employ locally supported basis functions, e.g., piecewise linear ``hat functions'', or splines. If $d\geq 2$, this approach requires the sample points to be selected in a structured way, for instance using sparse grids, see \cite{Barthelmann2000} and Section~\ref{sec:padapt}, or a mesh-based discretization of $\PP$ \cite{Ferranti2011}.
\item \textbf{Local unstructured.} A very simple, but nonetheless practical, way to enforce \eqref{eq:weights} is to construct the weights using a Voronoi tessellation \cite{Du2003} of $\PP$ based on the sample points $\PP_{train}$, i.e.,
\begin{equation}\label{eq:nearestneighbor}
\psi_\ell(\pp)=\begin{cases}
1\quad\text{if }\pp_\ell=\argmin_{\pp^\star\in\PP_{train}}\norm{\pp-\pp^\star},\\
0\quad\text{otherwise}.
\end{cases}
\end{equation}
This results in a nearest-neighbor approach, characterized by low accuracy, but also by a great flexibility. In fact, this strategy does not even require the poles to be matched. The main drawback of this approach is that it does not ``follow'' the evolution of the poles in $\pp$-space, resulting in limited predictive capabilities.
\end{itemize}

\begin{rmrk}\label{rem:passivity}
Depending on the application, the interpolation conditions \eqref{eq:weights} may be complemented by additional constraints to preserve important system properties, like stability, realness, or passivity \cite{Ferranti2011,Grivet2016}.
\end{rmrk}

\subsection{Parameter adaptivity}\label{sec:padapt}
Until now, we have assumed the parameter sample points $\PP_{train}$ to have been fixed in advance. However, in many situations, it proves extremely useful to have some kind of adaptivity included in the sampling of $\PP$, so that samples may be added only where the surrogate model is particularly inaccurate, e.g., in our case, near pole mismatches or where large interpolation errors occur. Still, it is quite difficult to devise adaptive strategies in non-intrusive MOR, especially if the number of parameters $d$ is large, since not much is known about the parametric dependence of the problem. In the context of sampling from high-dimensional parametric spaces, \emph{sparse grids} have been often employed in MOR as a way to alleviate the curse of dimensionality, see, e.g., \cite{Baur2011,Chen2015,Haasdonk2011}. Here, the focus is on adaptive sampling, and we propose a technique based on \emph{locally-refined sparse grids}, closely related to that considered in \cite{Alsayyari2019}, which, in turn, relies on some ideas from \cite{Ma2009,Pflger2010}. For simplicity, we carry out our construction in the case $\PP=[-1,1]^d$. Generalizations to more complicated parameter domains may be obtained by isomorphism.

Consider the nested ($\Gamma(n)\subseteq\Gamma(n+1)$ for all $n$) one-dimensional point sets
\begin{equation}\label{eq:sgone}
\Gamma(n)=\begin{cases}
\emptyset\quad&\text{if }n<0,\\
\{0\}\quad&\text{if }n=0,\\
\{2^{1-n}j\}_{j=-2^{n-1}}^{2^{n-1}}\quad&\text{if }n>0.
\end{cases}
\end{equation}
We extend this definition to multiple dimensions by tensorization: for any \emph{level index} $\nn=(n_1,\ldots,n_d)\in\xZ^d$, we define the corresponding \emph{tensor grid} $\Gamma(\nn)=\Gamma(n_1)\times\Gamma(n_2)\times\ldots\times\Gamma(n_d)$. It is useful to define the infinite point set
\begin{equation*}
{\bm\Phi}=\bigcup_{\nn\in\xZ^d}\Gamma(\nn)=\bigcup_{n=0}^\infty\Gamma(n)^d,
\end{equation*}
which is dense in $\PP$ (it coincides with the dyadic rationals in $\PP$) and also a superset of any tensor grid (by construction). We will choose the adaptive sampling points within $\bm\Phi$.

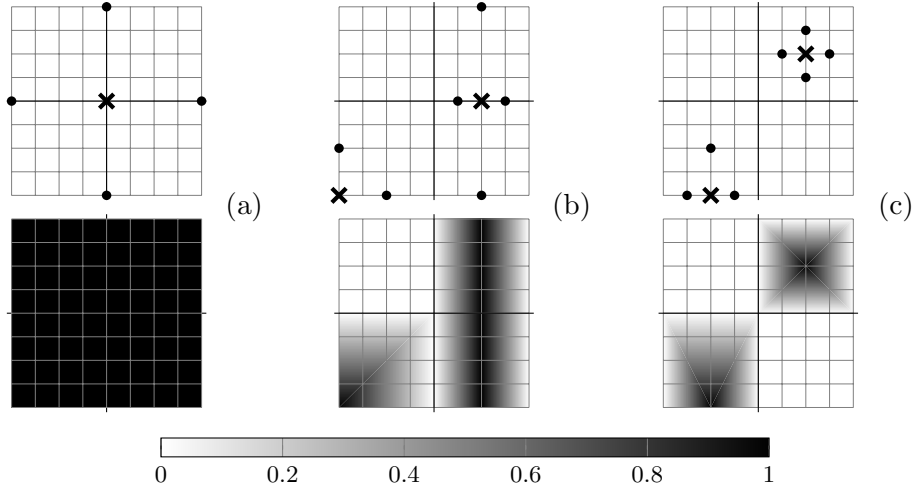
\begin{figure}[t!]
\centering
\includegraphics{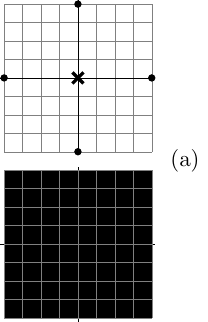}\hspace{.75cm}%
\includegraphics{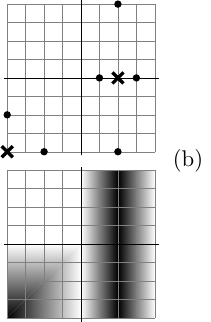}\hspace{.75cm}%
\includegraphics{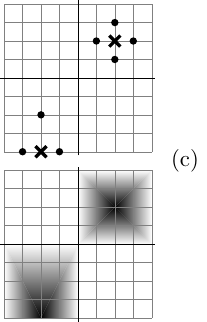}

\begin{tikzpicture}
\begin{axis}[
    hide axis, scale only axis,
    height=0pt, width=0pt,
	colorbar horizontal, colormap name = whiteblack,
    point meta min=0, point meta max=1,
    colorbar style={
        width=8cm, height =.25cm
    }
	]
    \addplot [draw=none] coordinates {(0,0)};
\end{axis}
\end{tikzpicture}
\caption{In the top row, the forward points of $(0,0)\in\Gamma(0,0)$ in (a), $(-1,-1)\in\Gamma(1,1)$ and $(1/2,0)\in\Gamma(2,0)$ in (b), and $(-1/2,-1)\in\Gamma(2,1)$ and $(1/2,1/2)\in\Gamma(2,2)$ in (c). In the bottom row, the corresponding hierarchical hat functions.}
\label{fig:sparsegrid}
\end{figure}

Now, assume that
\begin{equation}\label{eq:sparselevel}
\pp^\star\in\Gamma(\nn)\setminus\left(\bigcup_{k=1}^d\Gamma(n_1,\ldots,n_{k-1},n_k-1,n_{k+1},\ldots,n_d)\right),
\end{equation}
i.e., that the coordinates of $\pp^\star=\left(j_1/2^{n_1-1},\ldots,j_d/2^{n_d-1}\right)\in{\bm\Phi}$ are fractions in lowest terms (with $n_k=0$ if $j_k=0$). We define the \emph{forward points} of $\pp^\star$ as the ($\leq 2d$) elements of the discrete neighborhood
\begin{align*}
U(\pp^\star)=&\bigcup_{k=1}^d\left\{\pp\in\Gamma(n_1,\ldots,n_{k-1},n_k+1,n_{k+1},\ldots,n_d),\ \norm{\pp-\pp^\star}=2^{-n_k}\right\}\\
=&\PP\cap\bigcup_{k=1}^d\left\{\left(p_1^\star,\ldots,p_{k-1}^\star,p_k^\star\pm2^{-n_k},p_{k+1}^\star,\ldots,p_d^\star\right)\right\}.
\end{align*}
Moreover, to each $\pp^\star$ satisfying \eqref{eq:sparselevel}, we associate a \emph{hierarchical hat function} $\varphi_{\pp^\star}:\PP\to[0,1]$ according to the definition
\begin{equation}\label{eq:piecewiselinear}
\varphi_{\pp^\star}(\mathbf{x})=\prod_{k=1}^d\widehat{\varphi}_{p_k^\star,n_k}(x_k),
\end{equation}
with $\widehat{\varphi}_{p,0}(x)\equiv 1$ and, for $n=1,2,\ldots$,
\begin{equation*}
\widehat{\varphi}_{p,n}(x)=\begin{cases}
1-2^{n-1}\abs{x-p}\;&\text{if }\abs{x-p}<2^{1-n},\\
0\;&\text{if }\abs{x-p}\geq2^{1-n}.
\end{cases}
\end{equation*}
By construction, $\varphi_{\pp^\star}$ is zero at all $\pp'\in{\bm\Phi}$ of which $\pp^\star$ is a forward point (the \emph{backward points} of $\pp^\star$), and also at all backward points $\pp''\in{\bm\Phi}$ of such $\pp'$, etc., all the way back to $\pp=\bm 0$. We show some two-dimensional examples of forward points and of hierarchical hat functions in Figure~\ref{fig:sparsegrid}.

We rely on hierarchical hat functions to cast piecewise-linear interpolation problems over subsets of sparse grids. More precisely, given sample points $\PP^\star\subset{\bm\Phi}$, and data $\{f(\pp^\star)\}_{\pp^\star\in\PP^\star}$, the \emph{piecewise-linear interpolant of $f$ based on samples at $\PP^\star$} is the unique element $\widehat{f}_{\PP^\star}$ of $\text{span}\{\varphi_{\pp^\star}\}_{\pp^\star\in\PP^\star}$ which interpolates exactly the data: this means that there exist unique coefficients $\{c_{\pp^\star}\}_{\pp^\star\in\PP^\star}$, depending only on $\PP^\star$ and $\{f(\pp^\star)\}_{\pp^\star\in\PP^\star}$, such that
\begin{equation}\label{eq:piecewiseinterp}
f(\pp^\star)=\widehat{f}_{\PP^\star}(\pp^\star)=\sum_{\pp'\in\PP^\star}c_{\pp'}\varphi_{\pp'}(\pp^\star)\qquad\forall\pp^\star\in\PP^\star.
\end{equation}
The desired Lagrangian basis \eqref{eq:weights} can then be found by setting the data as $\{f(\pp^\star)=\delta_{\pp^\star\pp_l}\}_{\pp^\star\in\PP^\star}$, with $\delta$ the Kronecker delta. We note that the expression of each hierarchical basis function \eqref{eq:piecewiselinear} depends only on its support point $\pp^\star\in\PP^\star$, whereas the expression of each Lagrangian basis function \eqref{eq:weights} depends on the whole support set $\PP^\star$.

Now we are ready to describe our adaptive technique, which is summarized in Algorithm~\ref{algo:adapt}.
As in the greedy selection of frequency samples (Algorithm~\ref{algo:mri}), the adaptivity is achieved through a ``look-ahead'' idea, although here the approach is rather heuristic: we use the forward points of the current training set as test set, i.e., parameter values at which the accuracy of the current pROM is evaluated. If the surrogate model is too inaccurate at some of the test points, they are added to the training set. This loop is repeated until a specified tolerance is achieved at all current test points.

This approach, differently from the usual \emph{isotropic} adaptive sparse grid sampling \cite{Barthelmann2000,Nobile2008}, in general does not add whole levels $\Gamma(\nn)$, but only subsets of them. In fact, the training set is not even guaranteed to be \emph{downward-closed}, i.e., a point might be in the training set while some of its backward points are not. The matter of missing backward points is discussed to some detail in \cite[Section 3.2]{Alsayyari2019}. In the remainder of our presentation and in our numerical experiments, we do \emph{not} require missing backward points to be added to the training set, both for simplicity of exposition and (mainly) to reduce the cost of the offline phase\footnote{If one can afford a higher offline time, including backward points is advisable. However, we note that the increase in training cost could be significant, since each sparse grid point in $\xR^d$ has up to $2d$ backward points, and a (costly) frequency model must be built by MRI at each of them for error estimation.}. We remark that we are allowed to work with a non-downward-closed training set because the error estimator driving our $\pp$-adaptivity does not rely on interpreting the expansion coefficients $c_{\pp'}$ in \eqref{eq:piecewiseinterp} as ``hierarchical surpluses'', as is commonly done in adaptive sparse grids \cite{Ma2009}. As a side note, we observe that not including the backward points makes it necessary to recompute the expansion coefficients \eqref{eq:piecewiseinterp} from scratch whenever new training points are added.

Within each iteration, in order to quantify the accuracy of the pROM at a test parameter value $\pp$, we use the following strategy: 
\begin{enumerate}[(a)]
\item\label{item:greedy} We use the current pROM (whose training set does not include $\pp$, nor any of the other test points) to predict the frequency response at $\pp$, using \eqref{eq:ROMglobal}, where the weight functions $\psi_k$ are hierarchical hat functions \eqref{eq:piecewiselinear} built by piecewise-linear interpolation \eqref{eq:piecewiseinterp}, cf. Remark~\ref{rem:piecewise}.
\item Through Algorithm~\ref{algo:mri}, we build a frequency surrogate at $\pp$, which we take as ``truth frequency response'' at $\pp$. This requires solving the full model at $\pp$, at as many frequency points as required by the $z$-greedy procedure.
\item We compare poles and residues of the two models by employing \eqref{eq:matchfunc} as distance; this requires the solution of a pole-matching problem.
\end{enumerate}

For the sake of efficiency, it is crucial to observe that, over the different $\pp$-greedy iterations, function \textsc{BuildFrequencyROM} may be called multiple times with the same argument (not only when the pROM is built through \textsc{pROM\_Train}, see Algorithm~\ref{algo:build}, but also when evaluating the accuracy of the current model on the test set). As long as memory is not an issue, one should store frequency surrogates built at previous $\pp$-greedy steps, so that no expensive solve of the full model is wasted.

\begin{algorithm}[t]
    \caption{Adaptive parameter sampling}
    \label{algo:adapt}
    \begin{algorithmic}[1]
    		\Function{pROM\_Adapt}{}
    		\State Initialize $\PP_{train}\subset{\bm\Phi}$ (with ``few'' elements), tol $>0$
    		\Repeat \Comment{Greedy loop}
	    		\State pROM$\gets$ \Call{pROM\_Train}{$\PP_{train}$} \Comment{Algorithm~\ref{algo:build}}
	    		\State $\PP_{test}\gets\Call{ForwardPointsOf}{\PP_{train}}\setminus\PP_{train}$ \Comment{Update test set}
	    		\State $\PP_{train}^{new}\gets\{\ \}$
	    		\For{$\pp\in\PP_{test}$}
    				\State $\overline{\textup{ROM}}\gets\textup{pROM}(\cdot,\pp)$ \Comment{Predict surrogate at test point}
    				\State ROM$\gets$ \Call{BuildFrequencyROM}{$\pp$} \Comment{Build new surrogate at test point}
	    			\If{\Call{HeavisideDistance}{ROM, $\overline{\textup{ROM}}$} $>$ tol} \Comment{Distance defined in \eqref{eq:matchfunc}}
    					\State $\PP_{train}^{new}\gets\PP_{train}^{new}\cup\{\pp\}$
    				\EndIf
	    		\EndFor
	    		\State $\PP_{train}\gets\PP_{train}\cup\PP_{train}^{new}$
    		\Until{$\PP_{train}^{new}=\emptyset$}
    		\State \emph{Optional}:\ pROM$\gets$ \Call{pROM\_Train}{$\PP_{train}\cup\PP_{test}$} \Comment{Remark~\ref{rem:nowasteadapt}}
    		\State \textbf{return} pROM
    		\EndFunction
    \end{algorithmic}
\end{algorithm}

\begin{rmrk}\label{rem:piecewise}
In \eqref{item:greedy}, we have forced our pROM technique to reconstruct poles and residues only through piecewise-linear hat functions. However, as discussed towards the end of Section~\ref{sec:interp}, in some cases one may want to employ a matching-free nearest-neighbor reconstruction. This is easily achieved by using the piecewise constant basis \eqref{eq:nearestneighbor} instead of hierarchical hat functions. In this case, to better account for the approximation properties of interpolation basis, it is natural to employ Haar-type sparse grids \cite{Chkifa2014,Ma2009}, which can be obtained as in Section~\ref{sec:padapt}, replacing \eqref{eq:sgone} by
\begin{equation*}
\Gamma(n)=\begin{cases}
\emptyset\quad&\text{if }n\leq 0,\\
\{2^{1-n}j\}_{j=-2^{n-1}+1}^{2^{n-1}-1}\quad&\text{if }n>0.
\end{cases}
\end{equation*}
This essentially corresponds to restricting the sparse grid points to the \emph{interior} of $\PP$.
\end{rmrk}

Whatever the reconstruction strategy in step \eqref{item:greedy}, any of the methods presented in Section~\ref{sec:interp} can still be applied as a post-processing step, at the end of the greedy loop. The reason for this additional computation could be, for instance, a smoother representation of poles and residues, or the removal of eventual noise by regularization. To this aim, we wish to stress that some care should be used when selecting the pole-residue reconstruction strategy. Indeed, due to the local nature of the $\pp$-refinements, finer and coarser sampling regions may arise, which, if not taken into account, could lead to a poorly-behaving reconstruction. For instance, global polynomial interpolation over sampling points which are ``too wild'' can be an extremely ill-posed problem, due to a large Lebesgue constant \cite{Phillips2003}, whereas polynomial regression with low enough degree can be expected to behave more nicely.

\begin{rmrk}
The strategy that we presented is heuristic. In particular, it does not guarantee that, at the end of the greedy loop, the tolerance will be attained over the whole parameter domain, since we are using a relatively small (and sparse) test set to quantify the approximation error. Representing (``sketching'') the parameter domain by the test set can be justified only by assuming the resolution of the test set to be sufficiently fine. However, in practice, this is usually computationally unfeasible (especially if the number of parameters is large, due to the curse of dimensionality).
\end{rmrk}

\begin{rmrk}
In \cite{YueAdapt}, a somewhat similar $p$-adaptive approach was proposed, which, however, can be applied only to the single-parameter case. While the adaptivity there was essentially ``unidirectional'', adding samples progressively from one end of $\mathcal{P}$ to the other, here, through sparse grids, our approach acts more ``isotropically''.
\end{rmrk}

\begin{rmrk}\label{rem:nowasteadapt}
As in the frequency-adaptivity, see Remark~\ref{rem:nowaste}, once the greedy iterations are over, we can take advantage of the extra samples taken at test parameter points, and build a much richer pROM than the one that satisfied the tolerance constraint. Here, this idea is even more attractive than for MRI, since the test set can (and usually does) contain quite a large number of parameter values, as opposed to just 1.
\end{rmrk}

\section{Remarks on the smoothness of the Heaviside decomposition}\label{sec:bifurcation}
Based on how smoothly the system matrices $\{A,B,C,E\}$ in \eqref{eq:dynamical} depend on $\pp$, it is possible for the spectral quantities $\{\lambda^{(j)},m^{(j)},Y^{(j)}\}$ in \eqref{eq:output} to depend smoothly on $\pp$ as well. More precisely, continuous dependence is often passed on from matrices to Heaviside terms: small perturbations of the system matrices yield small perturbations of the poles $\lambda^{(j)}$ and of the residues $Y^{(j)}$, at least as long as multiplicities $m^{(j)}$ are independent of $\pp$ and poles do not cross. However, inheritance of analytic dependence cannot in general be guaranteed, since (polynomial) branches may naturally arise when poles cross. We refer to \cite[Chapter 12]{Reed1981} for an introductory discussion on the topic, and we report here two simple representative examples.

\subsection{A toy example of mode steering}\label{sec:ex1}
First, we showcase some of the intrinsic difficulties in dealing with crossing or almost crossing poles, even in the absence of bifurcations. The example below was obtained by generalizing a numerical test from \cite{Amsallem2011}.

For some fixed $\varepsilon\in\xR$, set $d=1$ and take
\begin{equation}\label{eq:ex1mat}
A_p=\begin{bmatrix}
2p & \varepsilon \\ \varepsilon & 0
\end{bmatrix},\quad B_p=\begin{bmatrix}
1 \\ 0
\end{bmatrix},\quad\text{and}\quad C_p=E_p=\begin{bmatrix}
1 & 0 \\ 0 & 1
\end{bmatrix}
\end{equation}
as the matrices defining a parametric dynamical system of the form \eqref{eq:dynamical}, which depend smoothly on $p\in\xR$ (and $\varepsilon$). The system output can be explicitly computed as
\begin{equation*}
Y(z,p;\varepsilon)=\begin{bmatrix}
z-2p & -\varepsilon \\ -\varepsilon & z
\end{bmatrix}^{-1}\begin{bmatrix}
1 \\ 0
\end{bmatrix}=\frac{1}{z^2-2pz-\varepsilon^2}\begin{bmatrix}
z \\ p_2
\end{bmatrix}=\frac{Y_p^{(1)}}{z-\lambda_p^{(1)}}+\frac{Y_p^{(2)}}{z-\lambda_p^{(2)}},
\end{equation*}
where the poles and residues are
\begin{equation}\label{eq:ex1polesres}
\lambda_p^{(1,2)}=p\pm\sqrt{p^2+\varepsilon^2}\quad\text{and}\quad Y_p^{(1,2)}=\frac12\begin{bmatrix}
1\pm\frac{p}{\sqrt{p^2+\varepsilon^2}}\\\pm\frac{\varepsilon}{\sqrt{p^2+\varepsilon^2}}
\end{bmatrix}.
\end{equation}

From \eqref{eq:ex1polesres}, it is not difficult to conclude that, if $\varepsilon=0$, the two poles coincide for $p=0$. However, for fixed $\varepsilon$, this degeneracy does not have a negative impact on the smoothness of the Heaviside decomposition\footnote{Interestingly, the \emph{joint} dependence on $p$ and $\varepsilon$ is non-smooth. Indeed, by comparing \eqref{eq:ex1polesreszero} and $$Y(z,0;\varepsilon)=\frac{1/2}{z-\varepsilon}\begin{bmatrix}
1 \\ 1
\end{bmatrix}+\frac{1/2}{z+\varepsilon}\begin{bmatrix}
1 \\ -1
\end{bmatrix},$$ we observe that the residues are discontinuous at $(p,\varepsilon)=(0,0)$. However, it is important to note that, while discontinuous, the residues stay uniformly bounded, since the matrix $A_p$ is real symmetric, hence diagonalizable, for all $(p,\varepsilon)\in\xR^2$.}, since the poles remain simple:
\begin{equation}\label{eq:ex1polesreszero}
Y(z,p;\varepsilon=0)=\frac{1}{z-2p}\begin{bmatrix}
1 \\ 0 
\end{bmatrix}+\frac{1}{z}\begin{bmatrix}
0 \\ 0
\end{bmatrix}.
\end{equation}

\begin{figure}[ht!]
\centering
\includegraphics{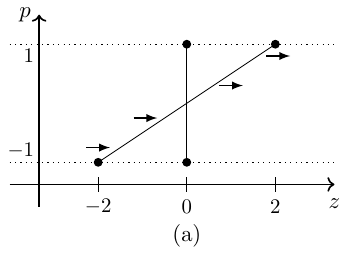}\hspace{.5cm}%
\includegraphics{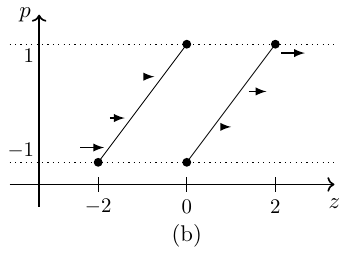}
\caption{Example of pole-matching for $\varepsilon=0$, with correct (a) and incorrect (b) matching; the residues at a few points, reconstructed by linear interpolation, are represented as small arrows. Missing arrows denote zero residues. Option (a) reconstructs the correct values exactly. Option (b) can happen only if $w=0$.}
\label{fig:ex11}
\end{figure}

We illustrate how the pole-matching algorithm performs in this simple example, using the exact Heaviside expansion \eqref{eq:ex1polesres} in place of the frequency surrogates one would obtain, e.g., via MRI. First, we fix $\varepsilon=0$, and consider the two parameter values $p_1=-1$ and $p_2=1$, where the Heaviside expansions of $Y$ are
\begin{equation}\label{eq:ex1polesresone}
Y(z,p_1;0)=\frac{1}{z+2}\begin{bmatrix}
1 \\ 0
\end{bmatrix}+\frac{1}{z}\begin{bmatrix}
0 \\ 0
\end{bmatrix}\quad\text{and}\quad Y(z,p_2;0)=\frac{1}{z-2}\begin{bmatrix}
1 \\ 0
\end{bmatrix}+\frac{1}{z}\begin{bmatrix}
0 \\ 0
\end{bmatrix}.
\end{equation}
We depict poles and residues in Figure~\ref{fig:ex11}. If the matching is accurate, pole $0$ should be matched with itself, and $-2$ with $2$.

As described in Section~\ref{sec:matching}, the matching criterion can be stated in terms of the cost matrix \eqref{eq:matchcost}, which here equals
\begin{equation}\label{eq:ex1matchcost}
D=\begin{bmatrix}
4 & 2+w\\2+w & 0
\end{bmatrix}.
\end{equation}
This means that matching $-2$ with $2$ and $0$ with itself has cost $4$, whereas matching $-2$ with $0$ and $0$ with $2$ has cost $4+2w$. Hence, as long as $w>0$, the algorithm performs the correct matching and recovers the exact response. However, if $w=0$, the two costs are the same, and any matching is allowed.

\begin{figure}[ht!]
\centering
\includegraphics{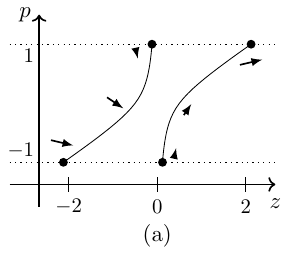}\hspace{.5cm}%
\includegraphics{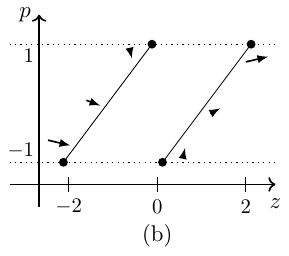}\hspace{.5cm}%
\includegraphics{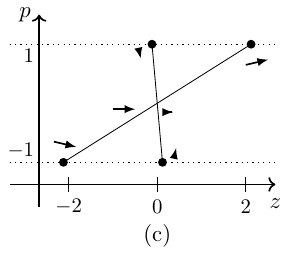}
\caption{Example of pole-matching for $\varepsilon=1/2$: exact solution (a), and pole-dominated (b) and residue-dominated (c) matching; the residues at a few points, reconstructed by linear interpolation for (b) and (c), are represented as small arrows.}
\label{fig:ex12}
\end{figure}

The matching becomes less trivial if $\varepsilon\neq 0$. For instance, let $\varepsilon=1/2$. We represent graphically poles and residues in Figure~\ref{fig:ex12} (a). By building the cost matrix \eqref{eq:matchcost} in this case (we omit the calculation here), we can conclude that the optimal matching changes depending on whether $w\mathrel{\ooalign{\raisebox{.6ex}{$>$}\cr\raisebox{-.6ex}{$<$}}}5-2\sqrt{5}\simeq0.53$: if poles have more importance than residues (i.e., $w$ is small), the surrogate poles do not cross, see Figure~\ref{fig:ex12} (b), whereas they do if the weight of residues is dominant (i.e., $w$ is large), see Figure~\ref{fig:ex12} (c).

At least qualitatively, case (b) appears to be a slightly better approximation of the poles and residues of the system. Still, neither of the surrogates (b) and (c) identifies the pole-residue behavior in a satisfactory way, and only adding more sample points (starting from $p_3=0$) will allow for a significant improvement in the quality of the approximation.

\subsection{A toy example of bifurcation}\label{sec:ex1b}
In this section we show a deceptively simple example of bifurcation, which, in practical applications, may arise due to unfavorable spectral properties of the problem (e.g., local non-diagonalizability).

\begin{figure}[ht!]
\centering
\begin{tikzpicture}[scale=.9]
\draw[->, thick] (-.65, -.75) -- (-.65, 2.5);
\draw[->, thick] (-1, -.575) -- (3.25, -.575) node[below=.1cm] {$p$};
\draw[semithick] (-.55, 0.) -- (-.75, 0.) node[left] {$-1$};
\draw[semithick] (-.55, 1.) -- (-.75, 1.) node[left] {$0$};
\draw[semithick] (-.55, 2.) -- (-.75, 2.) node[left] {$1$};
\draw[semithick] (0., -.475) -- (0., -.675) node[below] {$-1$};
\draw[semithick] (1.25, -.475) -- (1.25, -.675) node[below] {$0$};
\draw[semithick] (2.5, -.475) -- (2.5, -.675) node[below] {$1$};
\draw[thick,dashed,smooth]  (0,2.09544511501033) -- (0.125,2.04880884817015) -- (0.25,2) -- (0.375,1.94868329805051) -- (0.5,1.89442719099992) -- (0.625,1.83666002653408) -- (0.75,1.77459666924148) -- (0.875,1.70710678118655) -- (1,1.63245553203368) -- (1.125,1.54772255750517) -- (1.25,1.44721359549996) -- (1.3125,1.38729833462074) -- (1.375,1.31622776601684) -- (1.40625,1.27386127875258) -- (1.4375,1.22360679774998) -- (1.46875,1.15811388300842) -- (1.5,1) -- (1.46875,0.841886116991581) -- (1.4375,0.776393202250021) -- (1.40625,0.726138721247417) -- (1.375,0.683772233983162) -- (1.3125,0.612701665379258) -- (1.25,0.552786404500042) -- (1.125,0.452277442494834) -- (1,0.367544467966324) -- (0.875,0.292893218813452) -- (0.75,0.225403330758517) -- (0.625,0.163339973465924) -- (0.5,0.105572809000084) -- (0.375,0.051316701949486) -- (0.25,0) -- (0.125,-0.048808848170152) -- (0,-0.095445115010332);
\draw[thick,dashed] (1.5,1.) -- (2.5,1.);
\draw[thick,smooth]  (2.5,1.89442719099992) -- (2.375,1.83666002653408) -- (2.25,1.77459666924148) -- (2.125,1.70710678118655) -- (2,1.63245553203368) -- (1.875,1.54772255750517) -- (1.75,1.44721359549996) -- (1.6875,1.38729833462074) -- (1.625,1.31622776601684) -- (1.59375,1.27386127875258) -- (1.5625,1.22360679774998) -- (1.53125,1.15811388300842) -- (1.5,1) -- (1.53125,0.841886116991581) -- (1.5625,0.776393202250021) -- (1.59375,0.726138721247417) -- (1.625,0.683772233983162) -- (1.6875,0.612701665379258) -- (1.75,0.552786404500042) -- (1.875,0.452277442494834) -- (2,0.367544467966324) -- (2.125,0.292893218813452) -- (2.25,0.225403330758517) -- (2.375,0.163339973465924) -- (2.5,0.105572809000084);
\draw[thick] (1.5,1.) -- (0.,1.);
\draw (2.6, 2.15) node {\small $\Re(z)$};
\draw (.1, 2.35) node {\small $\Im(z)$};
\draw (1.25, -1.45) node {(a)};
\end{tikzpicture}\hspace{.25cm}%
\begin{tikzpicture}[scale=.9]
\draw[->, thick] (-.65, -.75) -- (-.65, 2.5);
\draw[->, thick] (-1, -.575) -- (3.25, -.575) node[below=.1cm] {$p$};
\draw[semithick] (-.55, 0.) -- (-.75, 0.) node[left] {$-1$};
\draw[semithick] (-.55, 1.) -- (-.75, 1.) node[left] {$0$};
\draw[semithick] (-.55, 2.) -- (-.75, 2.) node[left] {$1$};
\draw[semithick] (0., -.475) -- (0., -.675) node[below] {$-1$};
\draw[semithick] (1.25, -.475) -- (1.25, -.675) node[below] {$0$};
\draw[semithick] (2.5, -.475) -- (2.5, -.675) node[below] {$1$};
\draw[thick] (2.5,1.89442719099992) -- (0.,1.) -- (2.5,0.105572809000084);
\draw[thick,dashed] (0.,2.09544511501033) -- (2.5,1.) -- (0.,-0.095445115010332);
\fill (0.,1.) circle (.075);
\fill (2.5,1.) circle (.075);
\fill (0.,2.09544511501033) circle (.075);
\fill (0.,-0.095445115010332) circle (.075);
\fill (2.5,1.89442719099992) circle (.075);
\fill (2.5,0.105572809000084) circle (.075);
\draw (2.6, 2.2) node {\small $\Re(z)$};
\draw (.1, 2.4) node {\small $\Im(z)$};
\draw (1.25, -1.45) node {(b)};
\end{tikzpicture}\hspace{.25cm}%
\begin{tikzpicture}[scale=.9]
\draw[->, thick] (-.65, -.75) -- (-.65, 2.5);
\draw[->, thick] (-1, -.575) -- (3.25, -.575) node[below=.1cm] {$p$};
\draw[semithick] (-.55, 0.) -- (-.75, 0.) node[left] {$-1$};
\draw[semithick] (-.55, 1.) -- (-.75, 1.) node[left] {$0$};
\draw[semithick] (-.55, 2.) -- (-.75, 2.) node[left] {$1$};
\draw[semithick] (0., -.475) -- (0., -.675) node[below] {$-1$};
\draw[semithick] (1.25, -.475) -- (1.25, -.675) node[below] {$0$};
\draw[semithick] (2.5, -.475) -- (2.5, -.675) node[below] {$1$};
\draw[thick] (2.5,1.894427) -- (1.875,1.547723) -- (1.71875,1.41833) -- (1.5625,1.223607) -- (1.40625,1) -- (1.5625,0.776393) -- (1.71875,0.58167) -- (1.875,0.452277) -- (2.5,0.105573);
\draw[thick] (1.40625,1.) -- (0.,1.);
\draw[thick,dashed] (0,2.09545) -- (0.625,1.83666) -- (1.25,1.447214) -- (1.40625,1.273861) -- (1.5625,1) -- (1.40625,0.726139) -- (1.25,0.552786) -- (0.625,0.16334) -- (0,-0.09545);
\draw[thick,dashed] (1.5625,1.) -- (2.5,1.);
\fill (0,1) circle (.05);
\fill (0.625,1) circle (.05);
\fill (1.25,1) circle (.05);
\fill (1.40625,1) circle (.05);
\fill (1.5625,0.776393) circle (.05);
\fill (1.71875,0.58167) circle (.05);
\fill (1.875,0.452277) circle (.05);
\fill (2.5,0.105573) circle (.05);
\fill (1.5625,1.223607) circle (.05);
\fill (1.71875,1.41833) circle (.05);
\fill (1.875,1.547723) circle (.05);
\fill (2.5,1.894427) circle (.05);
\fill (0,2.09545) circle (.05);
\fill (0.625,1.83666) circle (.05);
\fill (1.25,1.447214) circle (.05);
\fill (1.40625,1.273861) circle (.05);
\fill (1.5625,1) circle (.05);
\fill (1.71875,1) circle (.05);
\fill (1.875,1) circle (.05);
\fill (2.5,1) circle (.05);
\fill (0,-0.09545) circle (.05);
\fill (0.625,0.16334) circle (.05);
\fill (1.25,0.552786) circle (.05);
\fill (1.40625,0.726139) circle (.05);
\draw (2.6, 2.2) node {\small $\Re(z)$};
\draw (.1, 2.4) node {\small $\Im(z)$};
\draw (1.25, -1.45) node {(c)};
\end{tikzpicture}
\caption{Real (full line) and imaginary (dashed line) parts of the exact roots of $z^2-p+\varepsilon$, for $\varepsilon=1/5$, with respect to $p\in[-1,1]$ (a). Surrogate roots from combinations of 2 (b) and 8 (c) $z$-surrogates.}
\label{fig:bifurcation}
\end{figure}
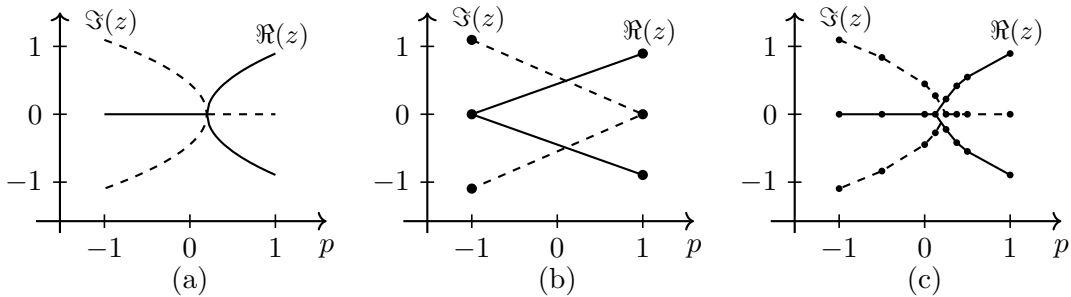

Fix a small $\varepsilon$, and consider the scalar problem with a single parameter
\begin{equation}\label{eq:ex0problem}
Y(z,p;\varepsilon)=\left(z^2-p+\varepsilon\right)^{-1}.
\end{equation}
For $p\neq\varepsilon$, the corresponding Heaviside expansion is readily found:
\begin{equation}\label{eq:ex0heaviside}
Y(z,p;\varepsilon)=\frac{1/(2\sqrt{p-\varepsilon})}{z-\sqrt{p-\varepsilon}}-\frac{1/(2\sqrt{p-\varepsilon})}{z+\sqrt{p-\varepsilon}},
\end{equation}
whereas, for $p=\varepsilon$, $Y(z,\varepsilon;\varepsilon)=z^{-2}$, with a double pole at 0. We can see that the poles are non-smooth (as functions of $p$) at $p=\varepsilon$, where they exhibit a bifurcation of degree 2, see Figure~\ref{fig:bifurcation} (a). Moreover, the residues are not only discontinuous at $p=\varepsilon$, but also unbounded there. This is due to the poles transitioning from single to double.

We apply our pMOR approach in this simple example, assuming, for simplicity, that the surrogate poles and residues coincide with the exact ones. If only the two parameter values $p_1=-1$ and $p_2=1$ are considered, the pROM will inevitably fail to identify $(p,z)=(\varepsilon,0)$ as a branch point. Indeed, our method ``follows'' separately the evolution with respect to $p$ of each pole. Consequently, the degree of each denominator in the Heaviside expansion \eqref{eq:ROMs} is kept equal to 1 even when crossing the singularity. The result, displayed in Figure~\ref{fig:bifurcation} (b), shows two surrogate pole lines which miss the branch point by ``twisting'' around it in the complex plane. Of course, as soon as more parameter points are added, the approximation quality improves significantly: in particular, if $p$-adaptivity is applied in this case, a progressively more accurate approximation of the exact poles is built as the greedy iterations proceed, see Figure~\ref{fig:bifurcation} (c).

It is important to remark that the structure of this kind of bifurcation would be much better identified by grouping together all of the involved branches, i.e., by merging several terms of the Heaviside decomposition into a single fraction with denominator degree $>1$ (here, 2) and coefficients smoothly dependent on the parameter, cf. the exact expression \eqref{eq:ex0problem}. This approach is implicitly applied by most methods based on global rational interpolation, e.g., \cite{Fevola2017,GTalocia2018,Ionita2014}, which, in fact, do not rely on the Heaviside decomposition of the output \eqref{eq:output}, but on its rational form
\begin{equation*}
Y(z,\pp)=C_\pp\left(zE_\pp-A_\pp\right)^{-1}B_\pp=\frac{P(z,\pp)}{Q(z,\pp)},
\end{equation*}
where $Q$ is a polynomial in $z$ of degree $\leq n_S$, with $\pp$-dependent coefficients. Procedures to integrate global rational interpolation concepts in our approach are being investigated, with the objective of improving the effectiveness in the approximation of non-smooth poles and residues. A possible idea revolves around a ``batch-matching'' of poles and residues (matching few-to-few instead of 1-to-1), which would allow an exact recovery of the simple bifurcation in the example above.

Before proceeding, we wish to briefly discuss what happens when a parameter sample point is added exactly at the branch point, which is the case, e.g., if $\varepsilon=0$ and we place a parameter sample at $p=0$. In such situations, there are two possibilities:
\begin{itemize}
\item The frequency surrogate (built, e.g., by MRI) identifies the double pole correctly. In this case, a simple Heaviside decomposition of the form \eqref{eq:ROMs} does not exist and our algorithm, as we presented it, fails. The batch-matching approach mentioned above (currently under investigation) would allow dealing in a natural way with this case.
\item The frequency surrogate mistakenly identifies the double pole as a couple of very close simple poles, whose residues have a large magnitude as a result, cf. \eqref{eq:ex0heaviside}. In this case, we may compute a Heaviside decomposition of the form \eqref{eq:ROMs}, and proceed as usual with the pMRI algorithm. However, due to the unboundedness of the residues, the approximation quality may be locally sub-optimal near the branch point.
\end{itemize}
We remark that, due to round-off noise (in the computation of the snapshots, in the MRI procedure, and in the Heaviside decomposition), the latter case is much more likely to present itself in practice.

\section{Numerical examples}\label{sec:examples}
We report in this section two numerical tests as evidence of the effectiveness of our technique. Our simulations were performed on the Helvetios cluster at EPFL \cite{Scitas}. For the sake of reproducibility, the corresponding code has been made available in \cite{Zenodo}.

\subsection{Laplacian eigenvalues on a parametric rectangle}\label{sec:ex2}
In this section we study a somewhat academic application in the field of PDEs for $d=1$, which was originally considered in \cite{Smetana2019}. Given $(k,H)\in\xC\times\xR^+$, we take the following Helmholtz equation on the rectangle $\Omega_H=\,]0,1[\,\times\,]0,H[$
\begin{equation}\label{eq:ex2helmholtzscaled}
\begin{cases}
-\left(\Delta+k^2\right)v_{k,H}(x_1,x_2)=f(x_1,x_2/H)\quad&\text{for }(x_1,x_2)\in\Omega_H,\\
v_{k,H}(x_1,0)=0\quad&\text{for }x_1\in\,]0,1[\\
\partial_{x_2} v_{k,H}(x_1,H)=\cos(\pi x_1)/H\quad&\text{for }x_1\in\,]0,1[\\
\partial_{x_1} v_{k,H}(0,x_2)=\partial_{x_1} v_{k,H}(1,x_2)=0\quad&\text{for }x_2\in\,]0,H[,
\end{cases}
\end{equation}
with $f:\,]0,1[^2\to\xR$ a piecewise constant forcing term (see \cite[Section 4.1]{Smetana2019} for the exact expression). Given how simply the geometry and the data of the problem depend on the parameter, we can recast \eqref{eq:ex2helmholtzscaled} on the reference domain $\Omega_1=\,]0,1[^2$ to make the parametric dependence emerge more clearly:
\begin{equation}\label{eq:ex2helmholtz}
\begin{cases}
-\left(\partial_{x_1x_1}^2+H^{-2}\partial_{x_2x_2}^2+k^2\right)u_{k,H}(x_1,x_2)=f(x_1,x_2)\quad&\text{for }(x_1,x_2)\in\Omega_1\\
u_{k,H}(x_1,0)=0\quad&\text{for }x_1\in\,]0,1[\\
\partial_{x_2} u_{k,H}(x_1,1)=\cos(\pi x_1)\quad&\text{for }x_1\in\,]0,1[\\
\partial_{x_1} u_{k,H}(0,x_2)=\partial_{x_1} u_{k,H}(1,x_2)=0\quad&\text{for }x_2\in\,]0,1[.
\end{cases}
\end{equation}
By inspection of the PDE above, we infer that $(z,p):=(k^2,H^{-2})$ is a good choice of parameters to study \eqref{eq:ex2helmholtz}. We set as frequency and parameter ranges $Z=[10,50]$ and $\mathcal{P}=[0.2,1.2]$, respectively.

After spatial discretization by FEM on a regular mesh with $n_S$ degrees of freedom, we obtain an algebraic problem of the form
\begin{equation}\label{eq:ex2helmholtzdiscrete}
\left(K_1+pK_2-zM\right)U(z,p)=b,
\end{equation}
where $K_1$, $K_2$, and $M$ are $n_S\times n_S$ matrices, and $U(z,p)$ and $b$ are vectors of size $n_S$ (here, we choose $n_S=10201$). We remark that \eqref{eq:ex2helmholtzdiscrete} is in the form of a dynamical system \eqref{eq:dynamical} with $B_p=b$. In order to conform to the functional setting of the PDE, we choose as norm over $V$ the functional $\xLtwo$ norm, which, in the discrete setting, corresponds to the energy norm induced by $M$: $\norm{\,\cdot\,}_V=\norm{M^{1/2}\,\cdot\,}_2$. The functional $\xHone_0$ or ($k$-weighted) $\xHone$ norms are also viable options; in our experience, the results of the simulation are barely affected by the choice of the norm.

\begin{figure}[p!]
\centering
\includegraphics[scale=.85]{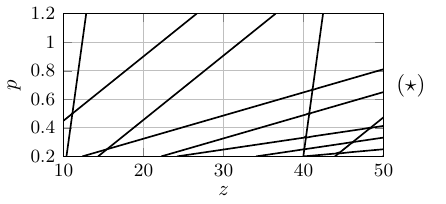}\hspace{.5cm}%
\includegraphics[scale=.85]{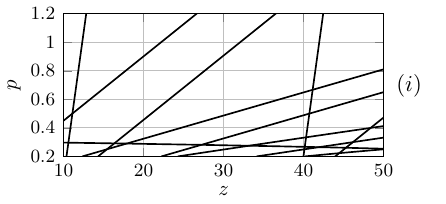}

\vspace{-.15cm}
\includegraphics[scale=.85, clip, trim = 0pt .25em 0pt 0pt]{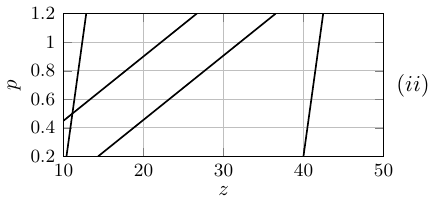}\hspace{.5cm}%
\includegraphics[scale=.85, clip, trim = 0pt .25em 0pt 0pt]{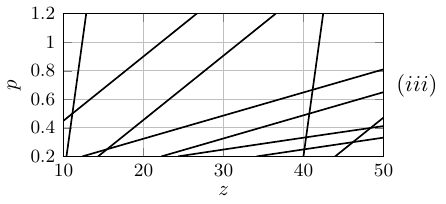}

\caption{Analytic poles \eqref{eq:ex2eigenvalues} of problem \eqref{eq:ex2helmholtzscaled} ($\star$). Surrogate poles obtained via double-greedy pMOR, with $\text{tol}_\text{synth}^{(i)}$, $\text{tol}_\text{synth}^{(ii)}$, and $\text{tol}_\text{synth}^{(iii)}$.}
\label{fig:ex21}
\end{figure}

We set as our target the FEM solution $U(z,p)$, i.e., we fix $C_p=I$, the identity matrix, in \eqref{eq:dynamical}. The poles of the Heaviside decomposition of $U$ gain additional importance as FEM approximations of the eigenvalues of the Laplace operator on the parametric domain $\Omega_H$. Due to our choice of domain, such eigenvalues\footnote{The set \eqref{eq:ex2eigenvalues} denotes the exact spectrum of the Laplace operator, without considering the FEM discretization. However, since the mesh is fine enough (the wavenumber $k$ is low and the mesh size is smaller than $k^{-2}$ \cite{Ihlenburg1995}), we expect the FEM eigenvalues and eigenvectors to be close to the analytic ones.} are actually available in closed form:
\begin{equation}\label{eq:ex2eigenvalues}
\left\{\pi^2k^2+\frac{\pi^2}{H^2}\left(\ell+\frac12\right)^2,\;k,\ell=0,1,\ldots\right\},
\end{equation}
see Figure~\ref{fig:ex21} ($\star$). We use the exact expression of the poles \eqref{eq:ex2eigenvalues} to validate the results obtained by our double-greedy pMOR technique.

\begin{figure}[p!]
\flushright
\includegraphics[scale=.85]{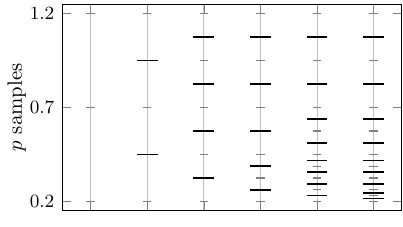}\hspace{-.21cm}%
\includegraphics[scale=.85]{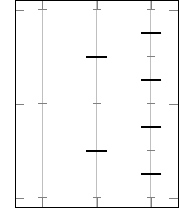}\hspace{-.21cm}%
\includegraphics[scale=.85]{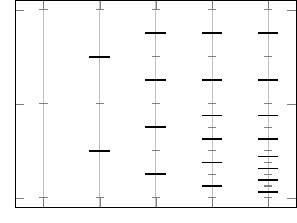}\hspace{2.5cm}

\vspace{-.325cm}
\includegraphics[scale=.85, clip, trim = 0pt .5em 0pt 0pt]{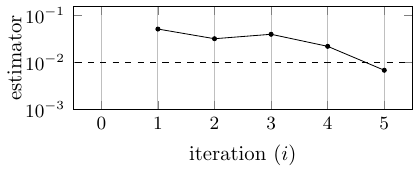}\hspace{-.21cm}%
\includegraphics[scale=.85, clip, trim = 0pt .5em 0pt 0pt]{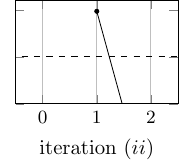}\hspace{-.21cm}%
\includegraphics[scale=.85, clip, trim = 0pt .5em 0pt 0pt]{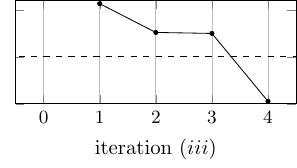}\hspace{2.5cm}

\centering
\caption{In the top row, new parameter samples added via greedy loop with $\text{tol}_\text{synth}^{(i)}$, $\text{tol}_\text{synth}^{(ii)}$, and $\text{tol}_\text{synth}^{(iii)}$. Thicker lines are used to denote the test set. The maximum error at each greedy iteration (measured by \textsc{HeavisideDistance} in Algorithm~\ref{algo:adapt}) is shown in the bottom row.}
\label{fig:ex22}
\end{figure}

\begin{figure}[p!]
\centering
\includegraphics[scale=.9, clip, trim = 0pt .25em 0pt -2pt]{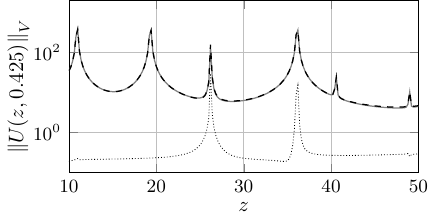}\hspace{.125cm}%
\includegraphics[scale=.9, clip, trim = 0pt .25em 0pt -2pt]{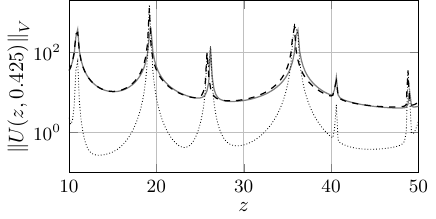}
\caption{On the left, the dashed line represents the surrogate solution norm $\|U\|$ obtained by double-greedy pMOR at $p^\star$. On the right, the dashed line represents the solution norm of the closest $z$-surrogate (at $p_\textrm{close}^\star$). In both cases, the solid line is the exact norm, whereas the dotted line is the absolute error $\|U-\textup{pROM}\|$.}
\label{fig:ex23}
\end{figure}

We employ the following computational setup:
\begin{itemize}
\item The frequency and parameter training sets are initialized as $3$ equispaced points in $Z$ and $\mathcal{P}$, respectively, whereas the frequency test set contains $100$ equispaced points in $Z$.
\item For MRI, Legendre polynomials are employed, whereas global monomials of degree 2 are used to interpolate (in a least squares sense) poles and residues after the greedy loop, in a post-processing ``compression'' step.
\item After both frequency and parameter greedy loops, the information at the test points is not wasted, but included in the final surrogates, see Remarks~\ref{rem:nowaste} and \ref{rem:nowasteadapt}. After the frequency greedy loops, we remove from the frequency surrogates any pole whose distance from $Z$ is larger than $20$, see Remark~\ref{rem:cleanup}.
\item The matching weight $w$ is set to $1$. The frequency and parameter greedy tolerances are set to $10^{-4}$ and $10^{-2}$, respectively.
\end{itemize}

On top of this, we consider three different choices for the tolerance $\text{tol}_\text{synth}$, which is employed to deal with unbalanced matching, see Algorithm~\ref{algo:unbalanced}:
\begin{itemize}
\item[$(i)$] $\text{tol}_\text{synth}^{(i)}=0$: all synthetic poles are kept.
\item[$(ii)$] $\text{tol}_\text{synth}^{(ii)}=1$: all synthetic poles are removed.
\item[$(iii)$] $\text{tol}_\text{synth}^{(iii)}=1/2$: synthetic poles are removed only if they account for the majority of the information about a pole.
\end{itemize}
Also, after the pole-matching has been completed, we improve the synthetic poles by extrapolating via global degree 2 monomials the non-synthetic poles, see Remark~\ref{rem:reconstruct}. This affects only $\text{tol}_\text{synth}^{(i,iii)}$.

We show the surrogate poles in Figure~\ref{fig:ex21}. In all cases, we can observe that the pole-crossings are handled well by the matching algorithm. This is likely due to the fact that residue information is taken into account ($w>0$), cf. Section~\ref{sec:ex1}.

In case ($i$), one erroneous pole crosses the frequency range $Z$ for small $p$. This is not caused by spurious poles being present in the frequency surrogates, but by an inaccurate matching of correct poles lying on different sides of $Z$. Due to the strict tolerance in case ($ii$), the pROM is blind to most of the poles which are not uniformly inside $Z$. Instead, the hybrid approach ($iii$) achieves a good compromise, missing only one pole that leaves $Z$ ``too quickly'' on the bottom right.

The ``hystory'' of Algorithm~\ref{algo:adapt}, namely the location of the new samples and the magnitude of the greedy error indicator, is portrayed in Figure~\ref{fig:ex22}. In the less strict cases ($i$) and ($iii$), the algorithm correctly identifies the ``busiest'' region of $\mathcal{P}$, i.e., small values of $p$, as critical for a good approximation. This results in local refinements near $p=0.2$. Instead, case ($ii$) only performs global refinements before terminating. This is actually a symptom of a general (undesirable) property of Algorithm~\ref{algo:adapt}: if $\text{tol}_\text{synth}$ is too large, pole-residue pairs might be removed too ``aggressively'' from the surrogate, so that the $\pp$-greedy termination criterion based on the Heaviside distance \eqref{eq:matchfunc} contains only a few terms. This, in turn, yields a smaller Heaviside distance, which is more likely to satisfy the prescribed $\pp$-greedy tolerance, potentially leading to an early termination. The simplest solution is to reduce $\text{tol}_\text{synth}$. Alternatively, one could partition the parametric domain $\PP=\bigcup_i\PP_i$, and then build a different pMOR surrogate on each parameter sub-domain $\PP_i$, see, e.g., \cite{Haasdonk2011}. On each sub-domain, it is less likely that a relevant pole will be discarded due to it being ``too often'' synthetic, cf. Section~\ref{sec:unbalanced}. This second approach is more costly (especially if $\PP$ is large and/or high-dimensional), but is particularly advantageous when the poles move very quickly through the parametric domain (i.e., if the gradient $\nabla_\pp\lambda$ is large).

In Figure~\ref{fig:ex23} (left), we compare exact and surrogate models at the point $p^\star=0.425$. For simplicity, we only consider the best surrogate, obtained with $\text{tol}_\text{synth}^{(iii)}$. The approximation seems of good quality, with pole locations and residue magnitudes being identified extremely well, and the approximation error is small. In particular, the results seem better than those obtained by using the closest frequency surrogate (i.e., the MRI built at the element of $\mathcal{P}_{train}$ closest to $p^\star$, namely, $p_\textrm{close}^\star=0.41875$), see Figure~\ref{fig:ex23} (right).

We also report a summary of the execution of the method and of the resulting pROM in Table~\ref{tab:ex2}. Our results agree with the main motivation behind the introduction of tol\textsubscript{synth}, namely that it should control how to deal with uncertain or missing information, resulting in richer (but also more noise-prone) or poorer (here, insufficient) surrogates.

\begin{table}[th!]
\small\centering
\begin{tabular}{|c|c|c|c|}
\cline{2-4}
\multicolumn{1}{c|}{} & ($i$) & ($ii$) & ($iii$) \\
\hline
tol\textsubscript{synth} & 0 & 1 & 1/2 \\
$p$-greedy iterations & 5 & 2 & 4 \\
$p$ samples & 19 & 9 & 17 \\
full model solves & 258 & 109 & 224 \\
number of poles (at each $p$) & 14 & 4 & 10 \\
synthetic poles (over all $z$-surrogates) & 81 (out of 266) & 0 (out of 36) & 22 (out of 170)\\
\hline
\end{tabular}
\vspace{.25em}\caption{Simulation details about the example in Section~\ref{sec:ex2}.}
\label{tab:ex2}
\end{table}

\subsection{Transmission line with high-dimensional parameter space}\label{sec:ex3}
\begin{figure}[ht!]
\centering
\includegraphics{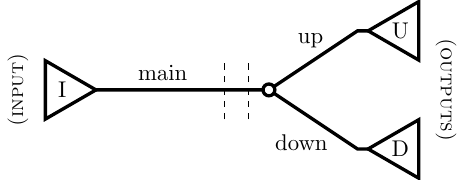}\hspace{1.25cm}%
\begin{circuitikz}[semithick]
\draw[dashed] (-3, -1.65) -- (-3, 1);
\draw[dashed] (1.75, -1.65) -- (1.75, 1);

\draw[white] (-3, -1.85) -- (1.75, -1.85);
\draw (-3, .25) to [R=$R$] (-1.125, .25) to [L=$L$] (.5, .25) -- (1.75, .25);
\draw (.5, .25) to [C=$C$] (.5, -1) node[ground]{};
\end{circuitikz}
\caption{Transmission line network diagram (left) and circuit representation of a unit cell (right).}
\label{fig:ex31}
\end{figure}
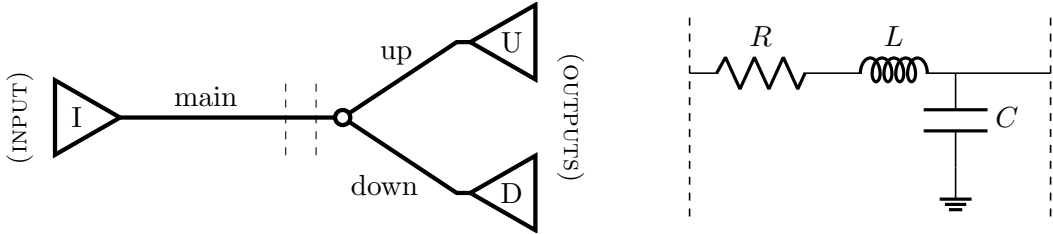

Our last numerical example concerns the analysis of the admittance parameters of the 3-port transmission tree depicted in Figure~\ref{fig:ex31}. A motivation and similar tests can be found, e.g., in \cite{Grivet2016}. Each branch is composed of a series of unit RLC cells: the ``main'' branch contains 400 cells, whereas the ``up'' and ``down'' branches contain 200 cells each. Resistance, inductance, and capacitance vary between cells: more precisely, if we restrict our focus to the main branch, the values of $R$, $L$, and $C$ of the $j$-th cell are, for all $j=1,\ldots,400$,
\begin{align*}
R^{\text{main},j}=&(1+\xi_R^{\text{main},j}+p_R^\text{main})\cdot 5\,\text{m}\Omega,\\
L^{\text{main},j}=&(1+\xi_L^{\text{main},j}+p_L^\text{main})\cdot 0.25\,\text{n}\textrm{H},\\
C^{\text{main},j}=&(1+\xi_C^{\text{main},j}+p_C^\text{main})\cdot 0.25\,\text{p}\textrm{F}.
\end{align*}
We employ $\xi$ to model random fluctuations of the nominal values in each cell, drawn from a uniform distribution with values between $-0.2$ and $0.2$. Such random values are fixed once and for all during the initialization of the model: we do not consider them as parameters in our analysis. Instead, the parameters $p$ denote branch-wide (independent of $j$) variations of the nominal values, and are envisioned to vary between $-0.1$ and $0.1$. For $j=1,\ldots,200$, the values of resistance, inductance, and capacitance in the secondary branches have the same expressions, with the subscript ``main'' being replaced by ``up'' and ``down'' for the top and bottom branches, respectively.

We consider the frequency range $Z=[0,8]$ GHz, and, as parameters in our pMOR approach, we take the 9-dimensional vector
\begin{equation*}
\pp=(p_R^\text{main},p_L^\text{main},p_C^\text{main},p_R^\text{up},p_L^\text{up},p_C^\text{up},p_R^\text{down},p_L^\text{down},p_C^\text{down})\in[-0.1,0.1]^9=\PP.
\end{equation*}
The admittance parameters can be found by solving a system of the form \eqref{eq:dynamical}, obtained by Modified Nodal Analysis: in particular,
\begin{align*}
A_\pp=&A^0+p_R^\text{main}A^\text{main}+p_R^\text{up}A^\text{up}+p_R^\text{down}A^\text{down},\\
E_\pp=&E^0+p_L^\text{main}E_L^\text{main}+p_L^\text{up}E_L^\text{up}+p_L^\text{down}E_L^\text{down}+p_C^\text{main}E_C^\text{main}+p_C^\text{up}E_C^\text{up}+p_C^\text{down}E_C^\text{down},
\end{align*}
while $B_\pp=B$ and $C_\pp=C$ are independent of $\pp$. In our case, the state $X$, which contains currents and voltages within the circuit, is a matrix of size $1603\times 3$, whereas the output $Y$ (the admittance matrix) has size $3\times 3$.

Our MOR setup for Algorithm~\ref{algo:adapt} is as follows:
\begin{itemize}
\item The frequency training set is initialized to the order $10$ Chebyshev points of $Z$, whereas the frequency test set contains $100$ equispaced points in $Z$; the parameter training set is initialized to $19$ points in $\PP$: the origin $\pp=\mathbf{0}$ and its 18 forward points.
\item For MRI, Legendre polynomials are employed, whereas piecewise linear hat functions are used to interpolate poles and residues.
\item After both frequency and parameter greedy loops, the information at the test points is not wasted, but included in the final surrogates, see Remarks~\ref{rem:nowaste} and \ref{rem:nowasteadapt}. After the frequency greedy loops, we remove from the frequency surrogates any pole whose distance from $Z$ is larger than $2$ GHz, see Remark~\ref{rem:cleanup}.
\item The matching weight $w$ is set to $1$; the frequency and parameter greedy tolerances are set to $10^{-3}$ and $10^{-2}$, respectively; the tolerance for unbalanced matching tol\textsubscript{synth} is set to 3/4.
\end{itemize}
We show a summary of the results of the offline training in Table~\ref{tab:ex3}.

\begin{figure}[p!]
\centering
\includegraphics[scale=.9, clip, trim = 0pt .25em 0pt 0pt]{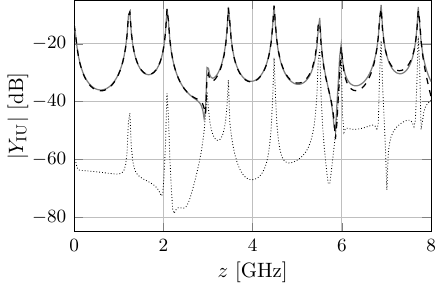}\hspace{.125cm}%
\includegraphics[scale=.9, clip, trim = 0pt .25em 0pt 0pt]{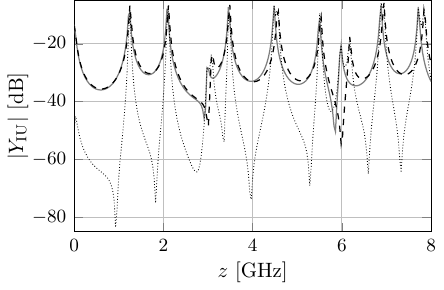}
\caption{On the left, the dashed line represents the surrogate magnitude of the admittance $Y_\text{IU}$ obtained by double-greedy pMOR at $\pp^\star$. On the right, the dashed line represents the admittance magnitude of the closest $z$-surrogate (at $\pp_\textrm{close}^\star$). In both cases, the exact magnitude is denoted by a solid line, and the absolute error by a dotted line.}
\label{fig:ex32}
\end{figure}

\begin{figure}[p!]
\centering
\includegraphics[scale=.95, clip, trim = 0pt .25em 0pt 0pt]{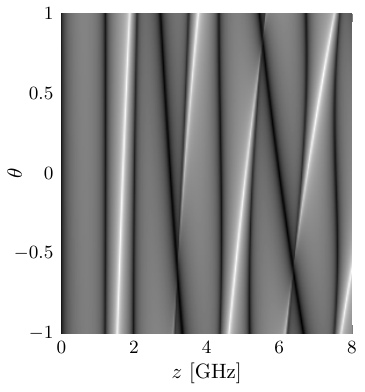}\hspace{.125cm}%
\includegraphics[scale=.95, clip, trim = 0pt .25em 0pt 0pt]{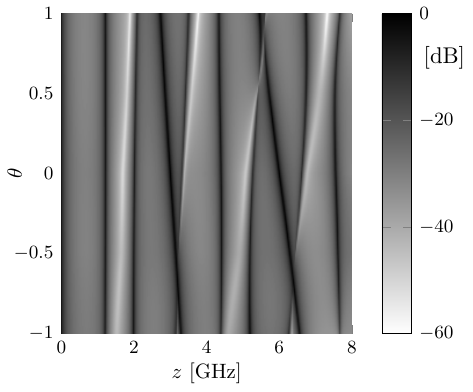}
\caption{Plot of $\abs{Y_\text{UD}}$ (in dB) for $\pp$ as in \eqref{eq:ex3line}: exact (left) and surrogate via pMOR (right). The color scale is the same for both plots.}
\label{fig:ex34}
\end{figure}

\begin{figure}[p!]
\centering
\includegraphics[scale=.95, clip, trim = 0pt .25em 0pt 0pt]{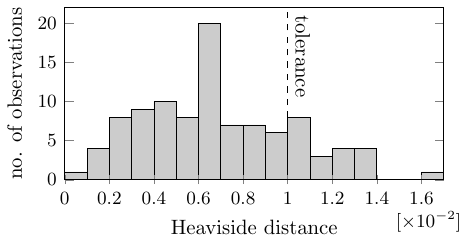}
\caption{Histogram of the $\pp$-greedy error indicator (measured by \textsc{HeavisideDistance} in Algorithm~\ref{algo:adapt}) between pMOR surrogate and local (high fidelity) MRI surrogate at 100 verification parameters, not belonging to the training or test sets of the pMOR surrogate.}
\label{fig:ex35}
\end{figure}

In Figure~\ref{fig:ex32} we compare exact and surrogate models at the randomly chosen point
\begin{equation}\label{eq:ex3rand}
\pp^\star=(-0.04, -0.08, -0.03, -0.01, 0.07, -0.09, -0.05, 0.03, 0.02).
\end{equation}
For simplicity, we only show the magnitude of the admittance between ports I and U. The quality of the approximation appears good, with pole locations and residue magnitudes being identified well. In particular, the results seem better than those obtained by using the closest frequency surrogate (i.e., the MRI built at the element of $\PP_{train}$ closest to $\pp^\star$, namely, $\pp_\textrm{close}^\star=[0,-0.1,0,0,0,-0.1,0,0,0]$).

We look at a more global picture in Figure~\ref{fig:ex34}. We let $\pp$ vary along a diagonal line between the two vertices $(0.1,0.1,0.1,-0.1,-0.1,\ldots,-0.1)$ and $(-0.1,-0.1,-0.1,0.1,0.1,\ldots,0.1)$:
\begin{equation}\label{eq:ex3line}
p_{R,L,C}^\text{main}=-0.1\theta\quad\text{and}\quad p_{R,L,C}^\text{up}=p_{R,L,C}^\text{down}=0.1\theta,
\end{equation}
with $\theta$ an auxiliary parameter. We plot the magnitude of the admittance between the two output ports for $z\in Z$ and $\theta\in[-1,1]$. Two different models are considered:
\begin{itemize}
\item The exact full order model \eqref{eq:dynamical}.
\item The surrogate obtained via double-greedy pMOR.
\end{itemize}

In Figure~\ref{fig:ex34}, we can observe that the exact admittance has 10 poles (darker lines) in the frequency range, some of which cross, and which, overall, create quite an intricate pattern. Still, the double-greedy pMOR surrogate seems to identify well the behavior of the quantity of interest, at least qualitatively. As could be expected, the quality of the approximation degrades slightly around pole intersections. A similar decrease in the accuracy of the surrogate can be observed also near the boundary of the parameter domain, since the two vertices of $\PP$ obtained for $\theta=-1$ and $\theta=1$ are \emph{not} elements of $\PP_{train}$.

To conclude the experiment, we also perform the following verification of the $\pp$-greedy tolerance employed to build the surrogate: we select 100 quasi-random (Halton) parameter points in $\PP$, outside the training and test sets of the pMOR surrogate. At each such point $\pp$, we build a reference frequency surrogate by $z$-adaptive MRI with the same parameters as above (10 Chebyshev points as starting training set, $z$-greedy tolerance of $10^{-3}$). Note that we only use \emph{new} snapshots at $\pp$ to build this reference model. Then we compare this model with the prediction given by our pMOR surrogate, obtained by simply plugging the value of $\pp$ in \eqref{eq:ROMglobal}. We make this comparison quantitative through the Heaviside distance \eqref{eq:matchfunc}, which is also the metric driving the $\pp$-adaptivity, cf. Algorithm~\ref{algo:adapt}. We plot the results in the form of a histogram in Figure~\ref{fig:ex35}. We can observe that 80\% of the verification points lie below the prescribed tolerance and that 99\% of them lie below 1.5 times the tolerance. Considering the discrete (sparse grid) nature of the test set, it is reasonable to expect that the tolerance will not be attained everywhere. In this context, we find our results satisfactory since they show that, even in the few cases where the tolerance is not satisfied, the error indicator is still within a small margin of the tolerance.

\begin{table}[th!]
\small\centering
\begin{tabular}{|c|c|c|}
\hline
$\pp$-greedy iterations & 6\\
$\pp$ samples & 457\\
full model solves & 9008\\
number of poles (at each $\pp$) & 14\\
synthetic poles (over all $z$-surrogates) & 183 (out of 6398)\\
\hline
\end{tabular}
\vspace{.25em}\caption{Simulation details about the example in Section~\ref{sec:ex3}.}
\label{tab:ex3}
\end{table}

\section{Conclusions and outlook}\label{sec:conclusions}
We have described the double-greedy pMOR approach for non-intrusive surrogate modeling of parametric problems, which samples adaptively in both frequency space and (potentially high-dimensional) parameter domain. In particular, the selection of parameter samples advances by trying to make the surrogate error small over a growing test set. We have illustrated with numerical examples the effectiveness of the method. Notably, we have shown that an accurate identification of number and behavior of poles and residues depends critically on some hyper-parameters ($w$ and tol\textsubscript{synth}) and, more generally, on the choice of a good strategy for interpolation over parameter space.

Among the several issues which remain unanswered we can find the following:
\begin{itemize}
\item Thanks to the degree of freedom provided by tol\textsubscript{synth}, the proposed heuristic strategy for unbalanced matching is reasonably flexible. Still, it is unclear whether an ``optimal'' choice of this tolerance exists and, if it does, how to find it for a given application, since it depends on quantities unavailable \emph{a priori}.
\item In Section~\ref{sec:bifurcation}, we have showcased some of the difficulties related to intersecting poles, in particular the potentially discontinuous behavior of residues. Applying the double-greedy approach in this case may yield inadequate surrogates, with the risk of a large number of iterations of the greedy loop. However, we have observed no such issues in our latter two numerical examples, despite multiple pole intersections. This is likely due to the beneficial spectral properties of the full order problems that we considered.
\item Issues similar to those discussed in the previous point are also possible (to a larger degree) in the case of multiple poles, even though, as discussed in Section~\ref{sec:bifurcation}, this case is quite unlikely to present itself thanks to numerical noise. On one hand, this problem should be partially addressed by MRI: for instance, when building the frequency surrogate, it should be possible to determine whether two poles are a noisy double pole or a couple of single poles. On the other hand, our pMOR algorithm should be extended to allow dealing with multiple poles, without compromising the overall complexity of the algorithm. Both of these directions are object of ongoing research.
\end{itemize}

\bibliographystyle{plain}
\bibliography{pMOR.bbl}
\end{document}